\documentclass[10pt]{article}
\usepackage{latexsym,amsmath,amscd,amssymb}
\makeatletter
\@addtoreset{figure}{section}
\def\thefigure{\thesection.\@arabic\c@figure}
\def\fps@figure{h,t}
\@addtoreset{table}{bsection}

\def\thetable{\thesection.\@arabic\c@table}
\def\fps@table{h, t}
\@addtoreset{equation}{section}

\makeatother

\begin{document}

\newtheorem{theorem}{Theorem}[section]
\newtheorem{definition}[theorem]{Definition}
\newtheorem{lemma}[theorem]{Lemma}
\newtheorem{remark}[theorem]{Remark}
\newtheorem{proposition}[theorem]{Proposition}
\newtheorem{corollary}[theorem]{Corollary}
\newtheorem{example}[theorem]{Example}
\newtheorem{examples}[theorem]{Examples}

\newcommand{\bfi}{\bfseries\itshape}

\newsavebox{\savepar}
\newenvironment{boxit}{\begin{lrbox}{\savepar}
\begin{minipage}[b]{15.8cm}}{\end{minipage}
\end{lrbox}\fbox{\usebox{\savepar}}}

\makeatletter
\title{{\bf Dynamics on Leibniz manifolds}}
\author{Juan--Pablo Ortega$^{1}$ and  V\'{\i}ctor
Planas--Bielsa$^{1}$}
\addtocounter{footnote}{1}
\footnotetext{Institut Nonlin\'eaire de Nice, UMR 129, CNRS-UNSA,
1361, route des Lucioles, 06560 Valbonne, France. }
\date{}
\makeatother
\maketitle

\addcontentsline{toc}{section}{Abstract}

\begin{abstract}
This paper shows that various relevant dynamical systems can be
described as vector fields associated to smooth functions via
a bracket that defines what we call a Leibniz structure. We show
that gradient flows, some dissipative systems, and
nonholonomically constrained simple mechanical systems, among
other dynamical behaviors, can be described using this
mathematical construction that generalizes the standard Poisson
bracket currently used in Hamiltonian mechanics. The symmetries
of these systems and the associated reduction procedures are
described in detail. A number of examples illustrate the
theoretical developments in the paper.
\end{abstract}


\section{Introduction}

It is well known that most of classical mechanics can be
formulated using a Poisson structure (see for
instance~\cite{abraham-marsden 1978, libermannmarle}, and
references therein). A Poisson structure on a manifold $P$ is a
bilinear map
$\{ \cdot ,
\cdot
\}: C^{\infty}(P)\times  C^{\infty}(P) \rightarrow 
C^{\infty}(P)$ that defines a Lie algebra structure on the algebra
$C^{\infty}(P)  $ of smooth functions on that manifold and that is 
a derivation on each entry. This property allows the association of
a vector field
$X _h
$ to any smooth function $h  \in C^{\infty} (P) $, usually referred
to as the Hamiltonian vector field of the Hamiltonian function
$h$. The properties of the bracket $\{ \cdot , \cdot\} $  have
important consequences on the dynamical features of the vector
field $X _h
$. For instance, its antisymmetry implies that the Hamiltonian
function
$h$  is a constant of the motion for $X _h
$. Additionally, the fact that $\{ \cdot , \cdot
\} $ satisfies the Jacobi identity implies that the flow of $X
_h  $ is a Poisson  map, that is, it respects the bracket. 

It has been noticed in recent times that a weakening of the
defining conditions of a Poisson system is sometimes necessary in
order to accommodate the description of more general dynamical
systems. A well known example is the use of brackets that do
not satisfy the Jacobi identity, known as almost Poisson brackets
(see Section~\ref{ejemplos}), in the context of nonholonomically
constrained mechanical systems.

In this paper we go all the way in this direction and we
work with a bracket, first introduced in~\cite{leibniz grabowski},  that
is just required to be bilinear and a derivation on each of its entries.
The derivation property, also known as the Leibniz rule, justifies why
we refer to this structure as {\bfi  Leibniz bracket}. The properties
exhibited by the  dynamical systems defined in this way are in general
very different from those presented by standard Hamiltonian systems since
most of the features of those systems are based on the Lie algebraic
properties of the bracket that we have chosen to drop. This construction
should not be mistaken with the Leibniz structures (also called
Loday algebras) introduced by Loday~\cite{loday} in the algebraic context.

The introduction  of the notion of Leibniz system, whose
elementary properties are presented in Section~\ref{Leibniz
systems section}, is justified by the great variety of relevant
systems whose natural underlying mathematical structure seems to
be based in this kind of brackets. A number of these systems are
described in Sections~\ref{ejemplos} and~\ref{Example: 
nonholonomic constraints}.  To be more specific, in
Section~\ref{ejemplos} we have identified the Leibniz structure
inherent to  a number of systems that can be found in the
literature and Section ~\ref{Example: 
nonholonomic constraints} is devoted to the Leibniz formulation
of simple mechanical systems subjected to nonholonomic
constraints. One of the main differences between our treatment
of this problem and other bracket formulations for nonholonomic
systems is the fact that our bracket is defined in the entire
phase space of the unconstrained system, unlike other approaches
that provide a bracket only on the constraint
submanifold. The Leibniz bracket that we construct in that
section associates to the Hamiltonian of the unconstrained system
a vector field that, when restricted to the constraint
submanifold, coincides with the evolution vector field of the
constrained system. As we point out in Section ~\ref{Example: 
nonholonomic constraints}, the construction of this bracket
involves the use of certain extensions that make it not to be
uniquely determined by the dynamics that we want to describe.
This leeway can be used in specific examples (see
Section~\ref{Example:  nonholonomic constraints}) to encode in a
single bracket entire families of constraints, which may be very
useful in  the study of bifurcation problems in the nonholonomic
context.

Section~\ref{Symmetries and reduction of Leibniz systems}
contains a first approach to the study of the symmetries and the
reduction of a Leibniz system. We introduce a notion of momentum
map associated to symmetries that do not necessarily preserve
the Leibniz bracket but that nevertheless produce non trivial
conservation laws as long as the Hamiltonian function is
invariant. We also formulate a theorem that spells out, under
certain regularity assumptions, how the orbit space of a symmetry
that does respect the Leibniz structure (in a weak sense that is
introduced in the text) is again a Leibniz manifold. This result
reproduces in the Leibniz context the well known result
for Poisson manifolds. The analog of symplectic or
Marsden--Weinstein reduction (see~\cite{mwr}) in this context is
the subject of ongoing research and will be treated elsewhere.

Finally, Section~\ref{The reduction of a presheaf of Leibniz
algebras}  contains the generalization to the Leibniz context of
the Poisson reduction results in~\cite{acg, poisson reduction
marsden ratiu, poisson reduction singular} that characterize the
situations in which a new Leibniz structure can be
obtained by restriction to a subset and projection to the orbit
space of a (pseudo)group of symmetries or to the leaf space of a
distribution.

\section{Leibniz systems}
\label{Leibniz systems section}

\begin{definition} 
\label{definition leibniz almost}
Let $P$ be a smooth manifold and let $C^{\infty}
(P)  $ be the ring of smooth functions on it. A {\bfi 
Leibniz bracket}  on  $P$  is a bilinear map $[ \cdot , \cdot ]:
C^{\infty}(P) \times C^{\infty}(P)\to
C^{\infty}(P)$ that is a derivation on each entry, that is,
\[
[fg,h]=[f,h]g+f[g,h]\text{\ \ \ and\ \ \ }[f,gh]=g[f,h]+h[f,g],
\]
for any $f,g,h \in C^{\infty}(P) $. We will say that the pair
$(P,[ \cdot , \cdot ])$ is a {\bfi 
Leibniz manifold}. If the bracket $[\cdot , \cdot ]$ is
antisymmetric, that is, it satisfies
$$
[f,g] =-[g,f] 
$$ 
for every pair of functions $f,g \in C^{\infty}(P)$ then we say
that $(P,[ \cdot , \cdot ])$ is an {\bfi 
almost Poisson manifold}. We will usually denote the almost
Poisson brackets with the symbol $\{ \cdot , \cdot \}  $.

A function $f  \in C^{\infty}(P) $ such that $[f,g]=0 $
(respectively, $[g,f]=0 $) for any $g \in C^{\infty} (P)  $ is
called a {\bfi  left} (respectively, {\bfi  right}) {\bfi 
Casimir} of the Leibniz manifold $(P, [ \cdot , \cdot ])$.
\end{definition}

\begin{definition} Let $(P,\{ \cdot , \cdot \})$ be an almost
Poisson manifold. We define  the {\bfi  Jacobiator} of the
bracket $\{ \cdot , \cdot \}  $ as the map
$
\mathfrak{J}:  C^{\infty}(P) \times  C^{\infty}(P) \times 
C^{\infty}(P) \to  C^{\infty}(P)
$ given by
\begin{equation}
\mathfrak{J}(f,g,h) = \{\{f , g\} ,h \}  +\{\{g ,h \} ,f
\}  +\{\{h ,f \} ,g \}. 
\label{Jacobiator}
\end{equation}
A {\bfi  Poisson structure} on $P$ is an
almost Poisson structure on $P$ for which the Jacobiator is the
zero map.
\end{definition}
The following lemmas are a direct consequence of the fact that
the Leibniz structure is a derivation.

Let $(P,[ \cdot , \cdot ])$ be a Leibniz manifold
and let $h$ be a smooth function on $P$. There exist two vector
fields $X^R_h $ and $X^L_h $ on $P$ uniquely characterized by the
relations
\[
X^R_h[f]=[f,h]\quad\text{and}\quad X^L_h[f](m)=-[h,f],
\quad\text{for any $f  \in  C^{\infty} (P) $.}
\]
Given two smooth functions $g,h \in C^{\infty}(P) $ there exists a
unique vector field
$X_{g,h}$ on
$P$ such that 
$$ X_{g,h}[f](m)=\mathfrak{J}(f,g,h)(m), \quad\text{for any $f
\in C^{\infty}(P)$}.
$$

We will call $X^R_h$ the {\bfi  Leibniz vector field}
associated to the {\bfi  Hamiltonian function} $h \in
C^{\infty}(P)$. In this paper, the abbreviation
$X_h$ will always denote
$X^R_h$. The flow $F _t $ of the vector field 
$X_h$ satisfies
\begin{equation}
\label{flow time evolution}
\left.\frac{d}{dt}\right|_{t=0}g(F _t (m))=[g,h](F _t(m)), \quad
\text{for any $g \in C^{\infty}(P)$}.
\end{equation}
A straightforward corollary of~(\ref{flow time evolution}) is
that in the context of almost Poisson manifolds the Hamiltonian
function is  a constant of motion, that is, if
$F _t$ is flow of $X _h  $ then $h \circ F _t  = h  $, for any $h
\in C^{\infty} (P)  $. Note that since
$[\cdot , \cdot]$ and $\mathfrak{J}$ are a derivation on each of
their arguments they only depend on the first derivatives of the
functions and thus, we can define two tensor maps
$B:T^{\ast}P \times T^{\ast}P \to \mathbb{R}$ and
$B_J:T ^{\ast}P \times T ^{\ast}P \times T ^{\ast}P \to
\mathbb{R}$  by 
$$ B(\mathbf{d}f,\mathbf{d}g) =[f,g] 
\qquad\text{and}\qquad
B_J(\mathbf{d}f,\mathbf{d}g,\mathbf{d}h)=\mathfrak{J}(f,g,h),
$$ 
for any $f,g,h \in C^{\infty}(P) $.
We can associate  to the tensor $B$ two vector
bundle maps $B_L^{\sharp}:  T^{\ast}P \to TP$ and $B_R^{\sharp} :
T^{\ast}P \to TP$ defined by the
relations
\[
B(\alpha, \beta)=- \langle \beta, B_L^{\sharp}(\alpha)\rangle
\quad \text{and}\quad B(\alpha, \beta)= \langle 
 \alpha, B_R^{\sharp}(\beta)\rangle 
\]
for any $\alpha, \beta \in T
^\ast  P $. Notice that when the bracket $[ \cdot , \cdot ] $ is
symmetric (respectively, antisymmetric) we have that
$B_R^{\sharp}=-B_L^{\sharp}
$  (respectively, $B_R^{\sharp}=B_L^{\sharp}
$) and
$X _h^R= -X _h ^L  $ (respectively, $X _h^R= X _h ^L  $), for any
$h \in C^{\infty} (P) 
$. We say that the Leibniz manifold $(P, [ \cdot
,\cdot ])$ is {\bfi  non degenerate} whenever the maps
$B_L^{\sharp} $ and $B_R^{\sharp} $  are vector bundle
isomorphisms.

\begin{definition} Let $(P, [ \cdot ,\cdot ])$ be a Leibniz
manifold. We define the {\bfi  left} and {\bfi  right
characteristic  distributions} $\mathcal{E}_L$ and
$\mathcal{E}_R$, respectively, by
$$
\mathcal{E}_L :={\rm span}\{ X_h^L  \, | \, h \in
C^{\infty}(P)\}=B_L^{\sharp}(T^{\ast}P)\quad\text{and}\quad
\mathcal{E}_R :={\rm span}\{ X_h^R  \, | \, h \in
C^{\infty}(P)\}=B_R^{\sharp}(T^{\ast}P).
$$
\end{definition}
Notice that if the Leibniz bracket $[\cdot , \cdot ]$ is
either symmetric or
antisymmetric then both distributions coincide. If additionally
the Leibniz manifold $(P, [ \cdot ,\cdot ])$ is non degenerate
then  $\mathcal{E}_L=\mathcal{E}_R=TP $  and 
we can define a  tensor field $\omega :TP \times TP \to
\mathbb{R}$ of type
$(0,2)$ on $P$ by
\begin{equation}
\omega (X_f , X_g )=[f , g], 
\label{leibniz form}
\end{equation} 
for any $f,g \in C^{\infty}(P) $. 
Given any point $m \in P $  and any vector subspace $V \subset T
_mP$ we denote
$$ V^{\omega}:=\lbrace w \in TP \mid 
\omega(m)(v,w)=0, \text{ for any } v \in V \rbrace.
$$
If the tensor $\omega $ is
antisymmetric (respectively, symmetric) then it is a two--form
(respectively, a pseudometric) on $P$. If additionally the form
$\omega $ is closed we say that $\omega $  is a {\bfi  symplectic
form} on $P$ and that the pair $(P, \omega)$ is a {\bfi 
symplectic manifold}.

Two functions smooth functions $h _1, h _2 \in C^{\infty}(P)$ 
on the Leibniz manifold $(P, [\cdot , \cdot ])$ are said to be
{\bfi  equivalent} if and only if
$[f,h_1 - h_2]=0
$, for any $f \in C^{\infty}(P) $ or equivalently, whenever $X
_{h _1}= X _{h _2} $. Notice that this definition establishes an
equivalence relation on the set
$C^{\infty}(P)$.

\begin{definition}
A {\bfi  Leibniz map}  between two Leibniz
manifolds
$(P_1, [\cdot , \cdot  ]_1)$ and $(P_2 , [ \cdot ,
\cdot ]_2)$ is a smooth map $\phi:P_1 \to P_2$ that satisfies
\[
\phi^{\ast}[f,g]_2=[\phi^{\ast}f,\phi^{\ast}g]_1, \quad \text{for
any $f , g
\in C^{\infty}(P_2)$.} 
\]
\end{definition}
\begin{lemma}
\label{equivariance} 
Let $\phi:(P_1, [ \cdot , \cdot]_1) \to (P_2, [
\cdot , \cdot]_2)$ be a Leibniz map. Let $h \in C^{\infty}(P_2)$
, $F^2_t$ be the flow of the Leibniz vector field $X_h$, 
$F^1_t$ the flow of
$X_{h \circ \phi}$, and ${\rm Dom}(F _t ^1)$ and ${\rm Dom}(F _t
^2)$ the domains of definition of $F _t ^1 $ and $F _t ^2 $,
respectively. Then ${\rm Dom}(F _t ^1)\subset \phi ^{-1} ({\rm
Dom}(F _t ^2))$ and 
\begin{equation}
\label{how flows go for leibniz}
F^2_t \circ \phi (z)= \phi \circ F^1_t (z),\text{ for any
$z$ in the domain  ${\rm Dom}(F _t ^1)$ of $F _t ^1$}.
\end{equation}
Additionally, $X_{h \circ \phi} $  and $X  _h  $ are $\phi
$--related, that is, $T\phi
\circ X_{h \circ \phi}=X_h \circ \phi$.
\end{lemma}

\noindent\textbf{Proof.\ \ } Let $z \in {\rm Dom}(F _t ^1)$ and $g
\in C^{\infty}(P_2)$ arbitrary. Using the Leibniz condition on
the map $\phi$  we can write
\begin{align*}
\mathbf{d}g ((\phi \circ F _t ^1)(z))&\cdot \frac{d}{dt} (\phi
\circ F _t ^1)(z) = \frac{d}{dt}g ((\phi \circ F _t ^1)(z))=
\frac{d}{dt}(g \circ  \phi)(F _t ^1 (z))=[g \circ \phi, h \circ 
\phi](F _t ^1 (z))\\
	&=[g , h ](\phi\circ F _t ^1) (z)= \mathbf{d} g ((\phi\circ F _t
^1) (z)) \cdot  X _h  ((\phi\circ F _t ^1) (z)).
\end{align*}
Since the function $g$  is arbitrary, this equality implies that
\[
\frac{d}{dt} (\phi
\circ F _t ^1)(z)=X _h  ((\phi\circ F _t ^1) (z)),
\]
which allows us to conclude
that $(\phi \circ F^1_t)(z)$ is an integral curve of $X_h$ through
the point $\phi(z )$. Since $F_t^2$ is the flow of $X _h  $
this automatically implies that $\phi(z) \in {\rm Dom}(F _t
^2 )$. As $z \in {\rm Dom}(F _t ^1 ) $ is arbitrary we get
that $\phi ({\rm Dom}(F _t
^1))\subset {\rm Dom}(F _t
^2 ) $ which implies that ${\rm Dom}(F _t ^1)\subset \phi ^{-1} ({\rm
Dom}(F _t ^2))$. Additionally, the uniqueness property of the
flow of a smooth vector field allows us to write that
$(\phi
\circ F^1_t)(z)=(F^2_t \circ \phi)(z)$. 

The $\phi$--relatedness of $X_{h \circ \phi} $  and $X  _h 
$ follows from taking the time derivative of~(\ref{how
flows go for leibniz}) at $t=0  $, recalling that ${\rm
Dom}(F _t ^1)$ becomes the entire manifold $P _1  $ when $t $
goes to zero.
\quad $\blacksquare$ 

\begin{proposition} 
\label{conserve categories}
Let $\phi:(P,\{\cdot , \cdot \}_P) \mapsto
(Q,[ \cdot , \cdot ]_Q)$ be a surjective Leibniz map. If
$(P,\{\cdot , \cdot \}_P)$ is a Poisson manifold then so is
$(Q,[ \cdot , \cdot ]_Q)$.
\end{proposition}

\noindent\textbf{Proof.\ \ }  Let $f,g
\in C^{\infty}(Q)$ arbitrary. The surjectivity of $\phi $
implies that any element in $Q$ can be written as $\phi(z) $,
for some $z \in P $. Hence,
\[
[f,g]_Q (\phi(z))=\{f \circ \phi, g \circ \phi\} _P (z)=-\{g \circ
\phi, f \circ \phi\} _P (z) =-[g,f]_Q (\phi(z)),
\]
which proves the antisymmetry of $[\cdot , \cdot ]_Q $.
Analogously, in order to prove that $[\cdot , \cdot ]_Q $
satisfies the Jacobi identity consider
$f,g,h
\in C^{\infty}(Q)$ and $z \in P $.  Since $\phi $ is a Leibniz
map we can write 
\begin{multline*}
[f,[g,h]_Q]_Q(\phi (z))+[g,[h,f]_Q]_Q(\phi
(z))+[h,[f,g]_Q]_Q (\phi
(z))\\
=\{f\circ \phi,\{g\circ
\phi,h\circ \phi\}_P\}_P+\{g\circ
\phi,\{h\circ \phi,f\circ \phi\}_P\}_P+\{h\circ
\phi,\{f\circ \phi,g\circ \phi\}_P\}_P (z)=0. \quad 
\blacksquare
\end{multline*}

\section{Examples}
\label{ejemplos}
\normalfont
As we already said symplectic and Poisson manifolds are
particular cases of Leibniz manifolds. We now briefly introduce
other non trivial  examples.

\smallskip

\noindent {\bf (i) Pseudometric brackets and gradient dynamical
systems.} Let $g:TP \times TP\rightarrow  \mathbb{R} $ be a
pseudometric on the smooth manifold $P$, that is, a symmetric non
degenerate tensor field of type $(0,2)$ on $P$. Let $g ^{\sharp}:
T ^\ast P  \rightarrow TP $ and $g ^{\flat}:T P\rightarrow T ^\ast
P $ be the associated vector bundle maps. Given any smooth
function $h \in C^{\infty} (P) $ we define its gradient $\nabla
h: P \rightarrow TP $ as the vector field on $P$ given by
$\nabla h:=g ^{\sharp} \mathbf{d} h$. Let $[\cdot , \cdot ]:
C^{\infty}(P) \times C^{\infty}(P) \rightarrow \mathbb{R} $ be
the Leibniz bracket defined by 
\[
[f,h]:= g (\nabla f, \nabla h),
\]
for any $f,h \in  C^{\infty}(P) $. We will refer to this bracket
as the {\bfi  pseudometric bracket} associated to $g $. This
bracket is clearly symmetric and non degenerate and the
Leibniz vector field $X _h  $ associated to any function $h
\in C^{\infty} (P) $ is such that  $X _h=\nabla h $.
These brackets are also called {\bfi Beltrami brackets}, 
see~\cite{crouch 81, van der Schaft-Beltrami}.
\smallskip

\noindent  {\bf (ii) The three--wave interaction.} A very
relevant problem in dynamics is the study of the interaction
between non linear oscillators and the energy exchange between
them. This problem can be viewed as an interaction between waves
of different frequencies with different resonance conditions. A
particular case that has deserved special attention is the so
called three--wave or triad interaction~\cite{alber luther
1998}. Following~\cite{Bloch} this problem can be formulated as
a dynamical system in
$\mathbb{R}^3$ that satisfies the differential equations given by
\[
\frac{d x}{dt}= s_1\gamma_1y z,\qquad \frac{dy}{dt}=
s_2\gamma_2 x z,\qquad \frac{dz}{dt}= s_3\gamma_3 xy,
\]
where the parameters $s _1, s _2, s _3 \in \{-1, 1\} $ and
$\gamma _1,
\gamma _2, \gamma _3 $ are real numbers that satisfy  $\gamma _1+
\gamma _2+ \gamma _3= 0 $.
This system happens to be a particular case of point {\bf (i)}
by taking the Leibniz bracket induced by the constant
pseudometric
$$ 
g=\begin{pmatrix} \frac{1}{s_1 \gamma_1} & 0 & 0\\ 0 &
\frac{-1}{s_2
\gamma_2} & 0\\ 0 & 0 & \frac{1}{s_3 \gamma_3}
\end{pmatrix}
$$
and the Hamiltonian function
$H(x,y,z)=xyz$.

\smallskip

\noindent {\bf (iii) Double bracket dissipation.} As we already
said the Leibniz dynamical systems induced by an almost
Poisson bracket are energy preserving. Brockett~\cite{brockett
1988, brockett 1993} has proposed the modelling of certain
dissipative phenomena by adding a symmetric bracket to a known
antisymmetric one, that is,
$$
[ \cdot , \cdot ]_{{\rm Leibniz}}=\{\cdot , \cdot\}_{{\rm skew}}
+[\cdot , \cdot]_{{\rm sym}}
$$
where the bracket $\{\cdot , \cdot\}_{{\rm skew}} $ is
skewsymmetric, $[\cdot , \cdot]_{{\rm sym}} $ is symmetric, and
hence the sum is a Leibniz bracket. This scheme allows the
modeling of a surprising number of physical examples. The reader
is encouraged to check with~\cite{lom, BlKrMaRa} for an account of
applications and references in this direction. 

A particularly simple example that fits into this framework is
the equation arising from the Landau--Lifschitz model for the
magnetization vector
${\bf M}$ in an external vector field ${\bf B}$,
\begin{equation}\label{dobr}
\dot{{\bf M}}=\gamma {\bf M} \times {\bf B} + \frac{\lambda}{\|{\bf M}\|^2} 
( {\bf M} \times (
{\bf M} \times {\bf B}))
\end{equation} where $\gamma$ and $\lambda$ are physical
parameters. This equation is Leibniz  in our sense if we take
the Leibniz bracket on $ \mathbb{R} ^3  $ given by the sum of the
two brackets
\[
\{f , g\}_{{\rm skew}}({\bf M})  :={\bf M} \cdot (\nabla f ({\bf M}) 
\times \nabla
g ({\bf M})) \quad \text{and} \quad
[f , g]_{{\rm sym}}({\bf M})  :=\frac{ \lambda ({\bf M} \times
\nabla f({\bf M}))({\bf M} \times
\nabla g({\bf M}))}{\gamma
\|{\bf M}\| ^2},
\]
where the symbol $\times  $ denotes the standard cross product on
$\mathbb{R}^3 $ and $\nabla $ is the Euclidean gradient. With this
bracket the differential equation~(\ref{dobr}) corresponds to the 
expression of the Leibniz vector field determined by the
function
\[
h ({\bf M})= \gamma {\bf B} \cdot  {\bf M}.
\]

Another related example is the differential equation satisfied by
a rigid body subjected to certain dissipation (see
~\cite{lom})
\begin{equation}
\label{equation solid dissipation}
\dot{{\bf M}}={\bf M} \times {\bf \Omega} + \alpha \left( {\bf M}
\times \left( {\bf M} \times {\bf \Omega} \right) \right).
\end{equation}
In this expression $ {\bf M} $ is the momentum vector of the
solid  and ${\bf  \Omega }$ its angular
velocity, both in body coordinates. Recall that 
\[{\bf \Omega}:=\left(\frac{M_1}{I_1},\frac{M_2}{I_2},
\frac{M_3}{I_3}\right),\]
where $(I _1, I _2,I_3)$  are the components of the inertia tensor
of the body with respect to a basis in which this tensor is
diagonal. If we take the same bracket as before with
$\alpha=\frac{\lambda}{\gamma \|{\bf M}\|^2}$, the
equations~(\ref{equation solid dissipation})  coincide with the
Hamilton equations corresponding to the function  
\[h({\bf M})=\frac{1}{2}
\left(\frac{M_1^2}{I_1}+\frac{M_2^2}{I_2}+\frac{M_3^2}{I_3} 
\right).\]

\smallskip 

\noindent {\bf (iv) Almost Poisson manifolds and non
holonomically constrained mechanical systems.} The equations of
motion of a simple mechanical system subjected to a constraint
can be written using D'Alembert's Principle. When the constraints
can be expressed as a linear function on the velocities these
equations admit a Leibniz formulation that corresponds to
the almost Poisson bracket introduced in
Definition~\ref{definition leibniz almost}. See~\cite{van der
Schaft and Maschke 1994, nonholonomic, cushman kemppainen
sniatycki, nonholonomic mechanical with bloch, almost Poisson,
sniatycki 2001, bloch-book}   and references therein. If the
constraints do not satisfy the linearity condition the almost
Poisson formulation ceases to be valid in general. Nevertheless,
if the constraints are affine on the velocities the problem still
admits a formulation in the context of Leibniz manifolds using a
bracket that in general is not antisymmetric. We discuss this
point in detail in the following section.

\section{Example:  nonholonomic constraints} 
\label{Example:  nonholonomic constraints} 

Let us consider a  simple mechanical
system characterized by a hyperregular Lagrangian
$L:TQ
\rightarrow
\mathbb{R}$ on the tangent bundle $TQ $ of a configuration space
$Q$. One way to impose kinematic constraints on that system
consists of fixing an affine subbundle  $C
\subset TQ$, usually referred to as the set of admissible
kinematical states. The
hyperregularity of $L$ implies that the associated Legendre transform
$\mathbb{F}L:TQ \rightarrow T ^\ast Q $ is a diffeomorphism
that can be used to define an associated Hamiltonian dynamical
system on $T ^\ast  Q $ with Hamiltonian function $H$, as well as
the {\bfi   Hamiltonian constraint submanifold} 
$D:=\mathbb{FL}(C) $ on $T ^\ast Q $.  D'Alembert's Principle
defines (see~\cite{nonholonomic}) a vector subbundle $W \subset T
_D (T ^\ast Q)$   such that
if $TD\cap W=\{0\} $  and the Hamiltonian vector field $X _H| _D 
$ is a section of $TD\oplus W $ then the corresponding splitting
$$ 
X_H \vert_D=X_D^H + X_W ^H
$$
is well defined and $X_D^H $ is the vector field whose flow
describes the motion of the constrained dynamical system. The
vector field $X_D^H $ is usually referred to as the {\bfi 
evolution vector field} and the complementary  vector field
$X _W ^H 
$  as the {\bfi   constraint force field}.

Our goal in the following paragraphs consists of endowing
$T^{\ast}Q$ with  a Leibniz
structure $[ \cdot , \cdot ]$ such that the Leibniz vector field $X^R_H$
associated
to $H$ is such that
$$
X^R_H(z)=X_D^H(z) \qquad \forall z \in D
$$

\begin{theorem}
Assume $T^{\ast}Q$ paracompact. Let $H \in C^{\infty}(T^{\ast}Q)$ be a
smooth 
function
with Hamiltonian vector field associated $X_H$. Let $D \subset T^{\ast}Q$
be a 
closed and embedded constraint submanifold and 
$W \subset T_D(T^{\ast}Q)$ a smooth vector subbundle such that
$T_D(T^{\ast}Q)=TD \oplus W$. There exists a Leibniz structure 
$[ \cdot , \cdot ]$ on $T^{\ast}Q$ such that
\begin{equation}
\label{extended bracket}
X_H^R(z)=\pi X_H(z)=:X_D^H(z) \quad z \in D
\end{equation}
where $\pi:TD \oplus W \longrightarrow TD$ is the natural projection
\end{theorem}

\noindent\textbf{Proof.\ \ } Consider the bilinear mapping 
$$
[ \cdot , \cdot ]_D : C^{\infty}(T^{\ast} Q) \times C^{\infty}(T^{\ast} Q) 
\longrightarrow C^{\infty}(D)
$$
defined by $[f , g]_D(z):=\langle \mathbf{d}f(z) , \pi B^{\sharp}(z)
\left( 
\mathbf{d}g(z) \right)\rangle$, for any $z \in D$, and where
$B^{\sharp}: T^{\ast}\left( T^{\ast}Q\right)\longrightarrow
T^{\ast}\left( T^{\ast}Q\right)$ 
is the vector bundle isomorphism induced by the canonical symplectic form
of 
$T^{\ast}Q$. This bracket has a smooth section $\widetilde{B}_D:D
\longrightarrow 
\mathcal{T}^2_0 \left( T^{\ast} Q\right)$ associated given by
$$
\widetilde{B}_D(z)\left( \alpha_z , \beta_z\right):=\langle
\alpha_z, \pi B^{\sharp}
(z) (\beta_z) \rangle,
$$
for any $\alpha_z , \beta_z \in T^{\ast}_z\left(T^{\ast} Q\right)$. By the smooth 
Tietze extension theorem (see for instance~\cite{mta}, Theorem
5.5.9)
$\widetilde{B}_D$ can 
be extended to a smooth section 
$\widetilde{B}:T^{\ast}Q \longrightarrow \mathcal{T}_0^2 (T^{\ast}Q)$.
For any $f , g \in C^{\infty}\left( T^{\ast}Q\right)$ we define

$$
[f , g]_L(m)=\widetilde{B}(m)\left( \mathbf{d}f(m) , \mathbf{d}g(m)\right)
$$
The bracket 
$[ \cdot , \cdot]_L:C^{\infty} \left( T^{\ast}Q\right)\times C^{\infty}
\left( T^{\ast}Q\right) \longrightarrow C^{\infty} \left(
T^{\ast}Q\right)$ 
endows $T^{\ast}Q$ with a Leibniz structure.
Finally, we show that~(\ref{extended bracket}) holds. Indeed, for
any $z \in D$,
\[
X_H^R(z)=\widetilde{B}^{\sharp}(z) \left( \mathbf{d}H(z)\right)=
\widetilde{B}_D(z) \left( \mathbf{d}H(z)\right)=\pi B^{\sharp}(z) 
\left( \mathbf{d}H(z)\right)=X_D^H(z). \qquad \blacksquare
\]

\begin{remark}
\label{extension remark}
\normalfont
The Leibniz structures on $T^{\ast}Q$ for which~(\ref{extended
bracket})  holds are not unique.
This freedom in the construction of the bracket can be used in specific 
applications to study families of systems instead of just a particular
one, 
which may be of relevance in bifurcation theoretical problems. We make
this 
comment more specific with the following elementary example.

Consider a particle of mass $m$ constrained to move in a   
rotating hoop of mass $M$ whose axis of rotation is parallel
to the gravity and contains a diameter of the hoop. Let $\phi$ be
the angle that  parametrizes the position of the bead in the hoop
and
$\psi$ the angle  that characterizes the position of the hoop. This
setup can be seen as  a simple mechanical system with configuration
space the two torus
$\mathbb{T}^2$ and subjected to the affine constraint
$\dot{\psi}-\omega=0$,
where $\omega \in \mathbb{R}$ is the constant angular speed of the 
hoop.

From the point of view of the formalism that we introduced above, 
for each value $\omega$, there is a Hamiltonian constraint submanifold, 
$D_{\omega}$ and a vector subbundle
$W={\rm span} \left\{ \frac{\partial}{\partial p_{\psi}} \right\}$ 
provided by D'Alembert's principle (see ~\cite{nonholonomic}) such that
$$
T \left( T^{\ast} \mathbb{T}^2\right)=TD_{\omega} \oplus W.
$$
As in the proof of the theorem, there exists for each $\omega$ a 
smooth section $\widetilde{B}_{D_{\omega}} \longrightarrow 
\mathcal{T}_0^2 \left( T^{\ast} \mathbb{T}^2\right)$. 
Given that the spaces $D_{\omega}$, $\omega \in \mathbb{R}$ 
form a foliation of $T^{\ast}\mathbb{T}^2$ we can define 
the section $\widetilde{B}: T^{\ast}\mathbb{T}^2 \longrightarrow 
\mathcal{T}^2_0 \left( 
T^{\ast}\mathbb{T}^2\right)$ by $\widetilde{B}(z) = 
\widetilde{B}_{D_{\omega}}(z)$, where $z \in D_{\omega}$.
The bracket $[ \cdot , \cdot ]_L$ induced by $\widetilde{B}$ on 
$T^{\ast}\mathbb{T}^2$ can be used to describe the system for
{\it any} value of the angular velocity of the hoop. 
On other words if $H \in C^{\infty}(\mathbb{T}^2)$ is the
Hamiltonian of the system then $X_H^R(z)=X_{D_{\omega}}^H(z)$ 
for any $z \in D_{\omega}$ and any $\omega \in  \mathbb{R}$.
\end{remark}

\section{Symmetries and reduction of Leibniz systems}
\label{Symmetries and reduction of Leibniz systems}

The symmetries of a dynamical system are in general very useful
to simplify its study. In the particular case of
symplectic and Poisson manifolds this idea has been
specifically implemented using a procedure that is generically
known as {\bfi   reduction} (see~\cite{mwr, poisson reduction
marsden ratiu, poisson reduction singular, hsr}, an references
therein). In the next two sections we will adapt some aspects of
the reduction theory of symplectic and Poisson systems to the
context of Leibniz manifolds.

In this section we will consider symmetries of Leibniz systems
that are encoded under the form of Lie group and Lie algebra
actions. Let $\mathfrak{g}$ be a Lie algebra and $P$ a 
smooth manifold. We recall that a  {\bfi right (left)
Lie algebra action}
of $\mathfrak{g}$ on $P$ is a Lie algebra (anti)homomorphism $\xi \in
\mathfrak{g}\longmapsto \xi_P \in \mathfrak{X}(P)$ such that the 
mapping 
$(m, \xi)\in P
\times \mathfrak{g}\longmapsto \xi_P (m) \in TP$ is smooth. The
symbol $\mathfrak{X} (P)  $ denotes the set of smooth vector
fields on $P$. We will denote by $C^{\infty} (P) ^{ \mathfrak{g}
}$ the set of $\mathfrak{g} $--invariant smooth functions on
$P$, that is,  $C^{\infty} (P) ^{ \mathfrak{g}
}:=\{ f  \in C^{\infty}(P)\mid \mathbf{d} f \cdot  \xi_P=0  
\mbox{ for any }   \xi \in  \mathfrak{g}\}$

\begin{definition} 
Let $(P,[ \cdot , \cdot ])$ be a Leibniz
manifold and $B\in  {\cal T}^2_0(P) $ the associated Leibniz
tensor. Let  $G$ be a Lie group (respectively, $\mathfrak{g}$
a Lie algebra) acting on
$P$. We say that this $G$--action (respectively,
$\mathfrak{g}$--action) is a {\bfi  weak symmetry} of
$(P,[
\cdot ,
\cdot ])$ whenever
the algebra $C^{\infty}(P)^G$ (respectively, $C^{\infty}(P)
^{\mathfrak{g}} $) of
$G$--invariant functions on
$P$ is closed under the Leibniz bracket. We say that $G$
(respectively, $\mathfrak{g} $) is a {\bfi  strong symmetry}  if 
$G$ acts on
$P$ by Leibniz maps (respectively, if $\pounds _{\xi_P} B =0$,
for any $\xi \in \mathfrak{g}  $). Such actions will be sometimes
referred to as  {\bfi  canonical}.
\end{definition}

\begin{definition}
\label{momentum map leibniz category}
Let $(P,[ \cdot , \cdot ])$ be a Leibniz
manifold and $\mathfrak{g}$ a Lie algebra acting on $P$. We
say that the $\mathfrak{g}$--action on $P$   admits a {\bfi 
momentum map} $\mathbf{J}:P \rightarrow \mathfrak{g}^\ast$
whenever, for any $\xi \in  \mathfrak{g}$, there exists a smooth
function $f _\xi \in C^{\infty} (P)  $ such that the component
$\mathbf{J} ^\xi:= \langle \mathbf{J}, \xi \rangle  $ is also
smooth and 
\[
X _{\mathbf{J} ^\xi} ^L = f  _\xi \xi_P.
\]
We will call the function $f _\xi $  the $\xi$--{\bfi 
integrating factor}.
\end{definition}

The definition that we just introduced allows the formulation of
a {\bfi  Noether's Theorem} in the Leibniz context. Notice that,
unlike the classical result,  only the
invariance of the Hamiltonian function is required in the
following statement and that there are no compatibility
assumptions between the action and the Leibniz bracket.

\begin{proposition}
\label{noether one more non holonomic}
Let $(P,[ \cdot , \cdot ])$ be a Leibniz
manifold and $\mathfrak{g}$ a Lie algebra acting on $P$. Assume
that the $\mathfrak{g}$--action on $P$   admits a
momentum map $\mathbf{J}:P \rightarrow \mathfrak{g}^\ast$. Then
the level sets of the momentum map are preserved by the flows of
the  Leibniz vector fields associated to any
$\mathfrak{g}$--invariant function on $P$.
\end{proposition}

\noindent\textbf{Proof.\ \ } Let $h \in C^{\infty}(P)
^{\mathfrak{g}} $  and $\xi\in   \mathfrak{g} $ be arbitrary. Let
$F _t  $ be the flow of the Leibniz vector field $X _h^R $. Then
\[
\frac{d}{dt} \mathbf{J}^\xi\circ  F _t=X _h ^R[\mathbf{J}
^\xi]=[\mathbf{J} ^\xi, h]= X _{\mathbf{J} ^\xi} ^L [h]= f _\xi
\xi_P[h]=f _\xi
\mathbf{d} h
\cdot  \xi _P=0. \quad \quad\blacksquare
\]

\begin{remark}
\normalfont
The introduction of the integrating factors in the definition of
the momentum map is not motivated by particular needs of the
Leibniz category. Indeed, as we show in the following example,
even in the symplectic or in the Poisson category this seems to
be the only way to associate non trivial conservation laws to non
canonical symmetries.

Let $\mathbb{R}^2 _{u} $ be the upper half plane considered as a
symplectic manifold with form $\omega= \mathbf{d}x\wedge
\mathbf{d}y $. Let $(\mathbb{R},+)$ act on $\mathbb{R}^2 _{u} $
by  $a \cdot (x,y):= (x, e ^a y) $, for any $(x,y)\in 
\mathbb{R}^2 _{u} $ and any $a \in \mathbb{R} $. This action is
clearly not canonical but it still admits a momentum map with a
non trivial integrating factor that leads to a conservation law
via Proposition~\ref{noether one more non holonomic}. Indeed,
for any $\xi\in \mathbb{R} $, the infinitesimal generator $\xi
_{\mathbb{R}^2_{u}} $ is given by $ \xi
_{\mathbb{R}^2_{u}} (x,y)=(0, e ^{\xi}y)$, $(x,y)\in 
\mathbb{R}^2_{u} $. This vector field is not Hamiltonian and
hence this action does not have a traditional momentum map
associated. However, there is a momentum map available in the
sense of Definition~\ref{momentum map leibniz category} since
the maps $\mathbf{J}:\mathbb{R}^2_{u} \rightarrow  \mathbb{R}  $
and  $f  _\xi \in C^{\infty}(\mathbb{R}^2_{u}) $ given by
$\mathbf{J} (x,y):= x $  and $f _\xi=- \left(\xi/ e ^\xi y
\right) $, $(x,y) \in \mathbb{R}^2_{u} $, $\xi \in \mathbb{R}  $,
are such that  
$X _{\mathbf{J}^\xi}= f _\xi \xi _{\mathbb{R}^2_{u}} $. The
conservation law associated by Proposition~\ref{noether one more
non holonomic} to the existence of this momentum map can be phrased
by saying that the Hamiltonian flows corresponding to Hamiltonian
functions that depend only on the
$x$ variable preserve the vertical lines in $\mathbb{R}^2_{u} $.
\end{remark}

The following result shows that the orbit space corresponding
to the weak symmetry $G$ of  a Leibniz manifold $(P,[ \cdot ,
\cdot ])$ is also a Leibniz
manifold provided that the group action satisfies enough
regularity assumptions to guarantee that the quotient $P/G $ is a
regular quotient manifold, that is, $P/G $  can be endowed with
a (unique) smooth structure that makes the projection $\pi: P
\rightarrow P/G  $  a surjective submersion. This observation is
consistent with the result of Bates and
Sniatycki~\cite{bates sniatycki} in the context of nonholonomic
mechanics that shows that symmetry reduction for these systems
preserves the form of the equations of motion.

\begin{theorem}
\label{main reduction theorem for leibniz}
Let $(P,[ \cdot , \cdot])$ be a
Leibniz manifold and let $G$ be a Lie group acting on $P$ in such
way that the orbit space $P/G $ is a regular quotient manifold
(this is the case when, for instance, the action is free and
proper). Assume that
$G$ is a weak symmetry of
$(P,[
\cdot ,
\cdot])$. Then

\begin{description}
\item [(i)] $(P/G,[ \cdot , \cdot ]_{P/G})$ is a Leibniz manifold
with bracket $[ \cdot , \cdot ]_{P/G} $ uniquely determined by the
expression
\begin{equation} [f,g]_{P/G}\circ  \pi=[\pi^{\ast}f,\pi^{\ast}g
],
\label{redbra}
\end{equation} 
for any $f, g  \in  C^{\infty} (P/G)  $ and where $\pi: P \to P/G$
is the projection.
\item [(ii)] The Leibniz structure induced by the bracket $[
\cdot ,
\cdot]_{P/G}$ on $P/G$ is the only one for which the projection
$\pi: P \to P/G$ is a Leibniz map.
\item [(iii)]  Let $h\in C^{\infty}(P)^G$ be a smooth
$G$--invariant function on $P$ and $h
^{P/G}
\in C^{\infty}(P/G)$ the function on the quotient uniquely
determined by the expression
$
h^{P/G} \circ  \pi= h$. Let $X_h$ and $X _{h^{P/G}}$  be the
corresponding Leibniz vector fields on $(P ,[ \cdot , \cdot
] )$ and $(P/G,[ \cdot , \cdot ]_{P/G})$, respectively, and $F _t 
$ and $F _t ^{P/G}$ the associated flows. Then ${\rm Dom}(F _t)
\subset \pi^{-1} ({\rm Dom}(F _t ^{P/G})) $ and 
\[
F _t ^{P/G} \circ  \pi (z)= \pi \circ F _t (z),
\]
for any $z \in {\rm Dom}(F _t ) $. 
The vector fields $X_h$ and $X_{h^{P/G}}$ are
$\pi$--related.
\end{description}
\end{theorem}

\noindent\textbf{Proof.\ \ } {\bf (i)} We first check
that~(\ref{redbra}) is a good definition for the bracket $[\cdot
, \cdot ]_{P/G}$. Let
$m , m'$ be two points in $P $ such that $\pi(m) = \pi(m')$. This
equality implies that there exists an element $g \in G $ such
that  $m' = g \cdot   m $. Let now $f, g  \in  C^{\infty}(P)^G 
$  arbitrary.  By definition
$ [f,g]_{P/G}(\pi(m'))=[f \circ
\pi , g \circ \pi ](m')$. Since by hypothesis $C^{\infty}(P)^G$
is closed under the bracket  and $f \circ
\pi$ and $ g \circ \pi $ are $G$--invariant then so is
$[f \circ
\pi , g \circ \pi ]$ and hence
\[
[f,g]_{P/G}(\pi(m'))=[f \circ
\pi , g \circ \pi ](m')=[f \circ
\pi , g \circ \pi ](g \cdot m)=[f \circ
\pi , g \circ \pi ](m)=[f,g]_{P/G}(\pi(m)),
\]
as required.
The bracket $[\cdot , \cdot ]_{P/G}$ is clearly bilinear and is a
derivation on its two arguments. Therefore, $(P/G,[ \cdot , \cdot
]_{P/G})$ is a Leibniz manifold.
{\bf (ii)}  It is a consequence of the fact that the
projection
$\pi$ is a surjective submersion. {\bf (iii)} is a consequence of
{\bf (ii)} and Lemma~\ref{equivariance}. \quad $\blacksquare$  

\medskip

We emphasize that the weak symmetry condition on the Leibniz
bracket $(P, [\cdot , \cdot ])$ and the $G$--invariance of the
Hamiltonian $h$ do not suffice to ensure the equivariance of the
associated Leibniz flow $F _t $ of $X _h $. In general this is
only true whenever the symmetry is strong.

\begin{remark}
\normalfont
The second part of the theorem shows that if the $G$--action is a
weak symmetry of a Leibniz manifold $(P, [\cdot , \cdot ])$ then
the quotient $P/G $ admits a unique Leibniz structure $[\cdot ,
\cdot ]_{P/G}$ with respect to which the projection $\pi:P
\rightarrow P/G $ is  a Leibniz map. The converse is also true.
Indeed, let $f , g  \in C^{\infty}(P)^G  $ and  $\bar{f},
\bar{g} \in C^{\infty}(P/G)$ be the unique smooth functions such
that  $\bar{f} \circ \pi= f $  and  $ \bar{g} \circ \pi= g $.
Then the bracket $[f,g]\in C^{\infty} (P) $ is such that  for
any $h \in G $  and any $z \in P $:
\[
[f,g](h \cdot z)=[\bar{f} \circ \pi,\bar{g} \circ \pi](h \cdot
z)=[\bar{f} ,\bar{g} ]_{P/G}(\pi(h \cdot z))=[\bar{f} ,\bar{g}
]_{P/G}(\pi(z))=[\bar{f} \circ \pi,\bar{g} \circ \pi](z)=[f,g](z),
\] 
which proves that $[f,g]\in C^{\infty} (P)^G $ and hence that
the $G$--action is a weak symmetry.
\end{remark}

\begin{remark} 
\normalfont
Proposition~\ref{conserve categories} guarantees that if the
Leibniz manifold
$(P, [\cdot , \cdot])$ in the previous theorem is actually
Poisson then so is the reduced manifold $(P/G, [\cdot , \cdot
]_{P/G})$.
\end{remark}

\begin{examples}
\normalfont
\noindent {\bf (i) Double bracket dissipation.} 
Consider the
systems with double bracket dissipation that we studied
in~\S~\ref{ejemplos}. This time we will restrict our discussion to
the points in  $\mathbb{R} ^3 $ that do not lie in the third
axis. Assume that the magnetic vector field ${\bf
B}$ is constant and equal to the vector $(0,0,1)$ or, in the case
of the rigid body subjected to a dissipative interaction,  
assume that the moments of inertia
$I_1$ and
$I_2$ are equal. In both cases, the rotations 
around the third axis, which is a free group action on the restricted phase space,
leave
invariant  the Hamiltonian functions and constitute a strong
symmetry for the Leibniz system which allows us to
apply Theorem~\ref{main reduction theorem for leibniz}. For a
concrete realization of the quotient Leibniz structure we use the
invariant polynomials
$\sigma_1=\frac{1}{2}\left( M_1^2 + M_2^2)
\right)$ and $\sigma_2=M_3$ of this
action. In these coordinates the Leibniz
tensor associated to the reduced Leibniz structure takes the form
\[ B=\frac{\gamma}{2 \sigma_1 + \sigma_3^2}
\begin{pmatrix} -2 \sigma_1 \sigma_2^2 & 2 \sigma_1 \sigma_2\\ 2
\sigma_1 \sigma_2 & -2 \sigma_1
\end{pmatrix}.
\] 
The reduced Hamiltonian functions are $h(\sigma_1,
\sigma_2)=\gamma \sigma _2
$ in the first case and $h=\frac{\sigma_1}{I _1} +
\frac{\sigma_2^2}{2 I _3}$ in the second.

\smallskip

\noindent {\bf (ii) Reduction of a Poisson system with a
non--canonical symmetry.} 
Consider the Poison dynamical system $(\mathbb{R}^3_{\ast},\{
\cdot ,
\cdot
\},H)$ on  $ {\Bbb R}^3_{\ast}:=\mathbb{R}
^3\setminus\{(x,y,0)\in  \mathbb{R}^3\}$, where the Poisson
bracket
$\{
\cdot , \cdot \} $ is determined by the Poisson tensor
$$ B(x,y,z)=\begin{pmatrix} 0 & x & y\\ -x & 0  & x\\ -y & -x & 0
\end{pmatrix},
$$ 
for any $(x,y,z) \in \mathbb{R}^3 _{\ast}  $ and $H=\frac{1}{2}
\left(x^2 + y^2\right)$.
Let $G:=(\mathbb{R},+)$ act on ${\Bbb R}^3$ by the map $\phi:
\mathbb{R} \times \mathbb{R}^3\rightarrow {\Bbb R}^3$ given by
$$ a \cdot(x,y,z)=(x,y,e^a z), \quad\text{for any $a  \in
\mathbb{R}  $ and any $(x,y,z) \in {\Bbb R}^3$}.
$$ 
This action is not canonical since, for
example, $\phi_a^{\ast}\{y,z\}=x \not=e ^a x= \{ \phi^{\ast}_ay ,
\phi^{\ast}_az \}$. However, notice that since the algebra of
$G$--invariant functions is made by functions depending on just 
the first two variables, that is, $C^{\infty}(P)^G=\{ f \in
C^{\infty}(P)\mid f\equiv f (x,y)\}$ then $C^{\infty}(P)^G$ is
closed under the Poisson bracket. Consequently, by
Theorem~\ref{main reduction theorem for leibniz},
$(\mathbb{R}^3_{ \ast},\{ \cdot ,
\cdot
\},H)$ can be reduced by this action. The reduced Poisson space
is 
$\mathbb{R}^2$ with the Poisson structure given by the reduced
Poisson tensor
$$ B_{{\Bbb R}^3_\ast/G}(x,y)=\begin{pmatrix} 0 & x\\ -x & 0
\end{pmatrix}
$$ and the reduced Hamiltonian is $h (x,y)=\frac{1}{2} \left(x^2 +
y^2
\right)$.

\smallskip

\noindent {\bf (iii) A symmetry of the three--wave
interaction.} The Hamiltonian $H(x,y,z)=xyz  $ of the
three--wave interaction that we presented in
Section~\ref{ejemplos} is invariant with respect to the action of
the product of the two multiplicative groups 
$G:=\mathbb{R}^+\times
\mathbb{R}^+$ by
$(\lambda_1,\lambda_2) \cdot(x,y,z):=(\lambda_1
x,\lambda_2 y,\lambda_3 z)$, where $(x,y,z) \in \mathbb{R}^3 $, $
(\lambda_1 ,\lambda_2)\in \mathbb{R}^+\times
\mathbb{R}^+$, and  $\lambda_3:=(\lambda_1
\cdot \lambda_2)^{-1}$. The infinitesimal generator for this
action is given by
$\xi_{\mathbb{R} ^3}(x,y,z)=(ax,by,-(a+b)z)$, for any
$\xi:=(a,b) \in \mathbb{R}^2$. 
Even though this action is not even a weak Leibniz symmetry we
can associate to it a momentum map $\mathbf{J}: \mathbb{R}^3
\rightarrow  \mathbb{R} ^2 $ given by 
$$
\mathbf{J}(x,y,z)=\left(\frac{1}{2}\left(\frac{x^2}{s_1
\gamma_1}-\frac{z^2}{s_3 \gamma_3}\right),
\frac{-1}{2}\left(\frac{y^2}{s_2 \gamma_2}+\frac{z^2}{s_3
\gamma_3}\right) \right).
$$
By Proposition~\ref{noether one more non holonomic}, the
components of  ${\bf J}$ are constants of the motion for the flow
of the Leibniz vector field $X _H $.

\smallskip

\noindent {\bf (iv) Nonholonomically constrained particle.}
Consider a free particle in $\mathbb{R}^3$. We will encode this
setup as a Hamiltonian dynamical system on the cotangent bundle
$T ^\ast  \mathbb{R}^3 $ endowed with its canonical symplectic
structure. We will denote by $B \in \Lambda^2(T ^\ast   (T ^\ast
 \mathbb{R} ^3))$ the Poisson tensor associated to this
symplectic form. The Hamiltonian function of this system is
$H(x,y,z,p_x,p_y,p_z)=\frac{1}{2}\left(p_x^2 + p_y^2 + p_z^2
\right)$. Suppose that the particle is forced to satisfy the
affine constraint $\dot{x}+y \dot{z}-a=0$, where
$a\in   \mathbb{R}$. In this particular case, the Hamiltonian
constrained submanifold is given by:
$$ D_a=\lbrace (x,y,z,p_x,p_y,p_z)\in  T ^\ast \mathbb{R} ^3\mid
 p_x +y p_z -a=0 \rbrace.
$$ 
Consequently,
$$
T_{(x,y,z,p_x,p_y,p_z)}D  _a={\rm span \,}\left\{
\frac{\partial}{\partial x},\frac{\partial}{\partial
y}+p_z\frac{\partial}{\partial p_x},\frac{\partial} {\partial
z},\frac{\partial}{\partial p_y},-y \frac{\partial}{\partial p_x
}+\frac{\partial}{\partial p_z}\right\}
$$ 
Using D'Alembert's principle (see~\cite{nonholonomic}) we
choose the subbundle $W_a \subset T _{D _a}(T ^\ast 
\mathbb{R}^3)$ given by 
$$
W_a(x,y,z,p_x,p_y,p_z):={\rm span}\left\{
\frac{\partial}{\partial p_x}+y
\frac{\partial}{\partial p_z} \right\}
$$ 
that satisfies the regularity
condition
$T _{D _a}(T ^\ast 
\mathbb{R}^3)=T D _a\oplus W _a$.
We now follow the scheme introduced in Section~\ref{Example: 
nonholonomic constraints}.   A straightforward computation shows
that the projection $\pi_a:T D _a\oplus W _a \rightarrow T D
_a$  and the composition $\widetilde{B}_{D_a}^{\sharp} :=\pi_a
\circ B^{\sharp}$ are given, using canonical coordinates, by the
matrices
\begin{equation}
\label{matrices we like}
\pi(m)=\begin{pmatrix} 1 & 0 & 0 & 0 & 0 & 0\\ 0
& 1 & 0 & 0 & 0 & 0\\ 0 & 0 & 1 & 0 & 0 & 0 \\ 0 &
\frac{-p_z}{1+y^2} & 0 &
\frac{y^2}{1+y^2} & 0 & \frac{-y}{1+y^2}\\ 0 & 0 & 0 & 0 & 1 & 0\\
0 & \frac{-y p_z}{1 + y^2} & 0 & \frac{-y}{1+y^2} & 0 &
\frac{1}{1+y^2}
\end{pmatrix},
\qquad \widetilde{B}_{D_a}^{\sharp} (m)=\begin{pmatrix} 0 & 0 & 0
& 1 & 0 & 0\\ 0 & 0 & 0 & 0 & 1 & 0\\ 0 & 0 & 0 & 0 & 0 & 1 \\
 \frac{-y^2}{1+y^2} & 0 & \frac{y}{1+y^2} & 0 & 
\frac{-p_z}{1+y^2} & 0\\ 0 & -1 & 0 & 0 & 0 & 0\\
 \frac{y}{1+y^2} & 0 & \frac{-1}{1+y^2}&0 & \frac{-y p_z}{1 +
y^2} & 0 
\end{pmatrix},
\end{equation}
where $m=(x,y,z,p_x,p_y,p_z) $. Following the ideas
introduced in the Remark~\ref{extension remark} and noticing that the 
expression of $\widetilde{B}_{D_a}^{\sharp}$ does not depend on
the parameter $a$, we can trivially extend $\widetilde{B}_{D_a} $ 
to a Leibniz tensor $\widetilde{B} \in {\cal T}^2_0(T ^\ast
\mathbb{R}^3)$  whose restriction  to any $D _a $, $ a \in
\mathbb{R}  $, coincides with $\widetilde{B}_{D_a}$. Using this
extension and the expressions in~(\ref{matrices we like}) we can
write the evolution vector field of the constrained system as
$$ X_D^H=\widetilde{B} ^{ \sharp}\mathbf{d}H,
$$ which in canonical coordinates reads
\[
\dot{x}=p _x, \quad\dot{y}=p _y, \quad\dot{z}=p _z, \quad\dot{p
_x}=\frac{-p_z p_y}{1+y^2}, \quad\dot{p _y}=0, \quad\dot{p
_z}=\frac{-y p_z p_y}{1+y^2}.
\]
Note that the constraint does not need to be included in the
set of equations since for a set of initial conditions satisfying
the constraint, the dynamics will preserve it automatically. From
the point of view of the Leibniz formulation of the problem this
remark can be phrased by saying that the function defining the
constraint is a left Casimir of the Leibniz system $(T ^\ast
\mathbb{R} ^3, \widetilde{B})  $.

\smallskip

\noindent \textbf{A momentum map.} The Hamiltonian $H$ is
symmetric with respect to  the lifted action of the
translations on the configuration space $\mathbb{R} ^3  $.
The infinitesimal generators associated to this action are given
by
$\xi_{T ^\ast  \mathbb{R} ^3}(x,y,z,p_x,p_y,p_z)=(\xi,0)$, for
any $\xi \in \mathbb{R} ^3  $. The lifted translations along the
$OY $--axis admit a momentum map $\mathbf{J}: T ^\ast  \mathbb{R}
^3  \rightarrow \mathbb{R}  $ with respect to the Leibniz
structure $(T ^\ast
\mathbb{R} ^3, \widetilde{B})  $, given by 
$$
\mathbf{J}(x,y,z,p_x,p_y,p_z)=p_y + \phi(p_x + y p_z)
$$ where $\phi$ is an arbitrary smooth real valued function. By
Proposition
\ref{noether one more non holonomic} $\mathbf{J}$ is preserved by
the integral curves of the evolution vector field $X _D ^H $.

\smallskip

\noindent \textbf{Reduction.} Consider now the group
$G:=(\mathbb{R}^2,+)$ acting on  $T ^\ast \mathbb{R} ^3 $
by lifting  the translations in the
coordinates
$x$ and $z$. This action is a weak symmetry of $(T ^\ast
\mathbb{R} ^3, \widetilde{B})  $ and thus we can apply
Theorem~\ref{main reduction theorem for leibniz}.  Let 
$(T ^\ast
\mathbb{R} ^3/G,[ \cdot , \cdot ]_{T ^\ast
\mathbb{R} ^3/G},h)$ be  the reduced
system.   $T ^\ast
\mathbb{R} ^3/G$ can be identified with $\mathbb{R}^4 $ since the
points in this orbit space correspond to the elements of the form
$(y,p_x,p_y,p_z)$. Using this identification the
reduced Hamiltonian can be written as
$h=\frac{1}{2}
\left(p_x^2 + p_y^2 + p_z^2 \right)$ and the reduced Leibniz
tensor $\widetilde{B} _1  $ is given by 
$$ \widetilde{B} _1^{\sharp}(y,p_x,p_y,p_z)=\begin{pmatrix} 0 & 0
& 1 & 0\\
\frac{y}{1 + y^2} & 0 & \frac{-p_z}{1+y^2}& 0\\ 0 & 0 & 0 & 0\\
\frac{-1}{1+y^2}& 0 & \frac{-y p_z}{1 + y^2} & 0
\end{pmatrix}.
$$ 
This Leibniz structure has two independent left Casimirs
$C^l_1(y,p_x,p_y,p_z)=p_y \, , C^l_2(y,p_x,p_y,p_z)=p_x + y p_z$
and two independent right Casimirs
$C^r_1(y,p_x,p_y,p_z)=p_x \, , C^r_2(y,p_x,p_y,p_z)=p_z$.
Consequently the new Hamiltonian function
$\overline{h}:=h-\frac{1}{2} \left((C^r_1) ^2+  (C^r_2)
^2\right)$, that is,
$\overline{h}(y,p_x,p_y,p_z)=\frac{1}{2}p_y^2$, has the same
evolution vector field than that of $h  $. 
This
equivalent Hamiltonian admits the symmetry of translations in the
coordinates $p_x$ and $p _z $ and hence we can further reduce the
system onto a two dimensional Leibniz one  with tensor
$$
\widetilde{B} _2^{\sharp}(y,p_y)=\begin{pmatrix} 0 & 1\\
0 & 0
\end{pmatrix}
$$
and Hamiltonian function
$h_2=\frac{1}{2}p_y^2$.

\end{examples}

\section{The reduction of a presheaf of Leibniz algebras} 
\label{The reduction of a presheaf of Leibniz algebras} 

The reduction theorem that we presented in the previous section
contains extremely strong regularity hypotheses that allowed us to
have a smooth orbit space onto which the Leibniz bracket and the
corresponding equivariant dynamics can be dropped. When these
hypotheses are not present, the orbit space is not smooth anymore
but nevertheless, the Leibniz algebra, or more specifically, the
{\bfi  presheaf of Leibniz algebras}  associated to the bracket 
admits, under  certain circumstances, a projection to the
quotient. The algebraic  approach to reduction that we introduce
in the following paragraphs has its origins in the
works~\cite{acg, poisson reduction singular} carried out in the
context of symmetric Poisson manifolds.

We
recall that a {\bfi  sheaf} 
$ \mathcal{F}$ of functions on a topological space $P$ is a map that assigns
to any open
set $U$ a set of real valued functions $\mathcal{F}(U) $ which is an algebra 
under
multiplication. In the definition it is also required that for every inclusion
$V \subset U
$ of open sets there is a given homomorphism $ \mbox{\rm res} _V^U: \mathcal{F}(U)
\rightarrow  \mathcal{F}(V) $ called the {\bfi  restriction}
from $U$ to $V$ that satisfies the following conditions:
\begin{description}
\item [(SH1)] $\mathcal{F}(\emptyset)=\{0\} $ and $\mbox{\rm res} _U^U :
\mathcal{F}(U)
\rightarrow  \mathcal{F}(U)$ is the identity
map.
\item [(SH2)] If $W \subset V \subset U $ are open sets, then 
$\mbox{\rm res} _W^V \circ
\mbox{\rm res}_V^U= \mbox{\rm res} _W^U $.
\item [(SH3)] Let $U$ be an open set and $\{ V _i\}_{i \in I} $ 
an  open covering of
$U$. If $f \in  \mathcal{F}(U)$ is such that the restriction 
$\mbox{\rm res}^U_{V _i}(f)
$ of
$f$ to each $V _i $ is $0$, then $f=0 $
\item [(SH4)] Let $U$ be an open set, $\{ V _i\}_{i \in I} $  an
open  covering of
$U$, and let $f _i \in  \mathcal{F}(V _i) $ be given for each $i \in I $. 
Suppose that 
the restrictions of  $f _i $ and $f _j $ to $V _i\cap V _j$ are equal for 
all $i, j \in
I$. Then there exists a unique $f  \in  \mathcal{F}(U) $ whose restriction 
to each $V _i
$ is $f _i$ for all $i \in I $.
\end{description}
When the map $\mathcal{F} $ satisfies only properties {\bf (SH1)} and {\bf (SH2)}
we say
that $\mathcal{F} $ is a {\bfi  presheaf}. 
\index{presheaf}%
The elements in $ \mathcal{F}(U) $ are called
the {\bfi  sections}
\index{section!of a sheaf}%
\index{sheaf!section}%
of $\mathcal{F} $ over $U$. The elements in $\mathcal{F} (P)$ are called
{\bfi  global sections}.

\begin{definition}
\label{presheaf of Leibniz algebras definition} 
Let $M$ be a
topological space with a presheaf ${\cal F}$ of smooth functions.
A {\bfi  presheaf of Leibniz algebras} 
on $(P, \mathcal{F}) $ is a map $[ \cdot , \cdot ] $ that assigns
to each open set $U \subset M $ a bilinear operation $[ \cdot ,
\cdot ] _U :
\mathcal{F} (U) \times \mathcal{F}(U) \rightarrow \mathcal{F}(U)$
such that the pair
$(\mathcal{F}(U),[ \cdot , \cdot ] _U)$ is a  Leibniz algebra. A
presheaf of Leibniz algebras will be usually denoted as a triple
$(P, \mathcal{F}, [ \cdot , \cdot ])$.

We say that the presheaf of Leibniz algebras $(P, \mathcal{F}, [
\cdot , \cdot ])$ is {\bfi  non degenerate}  when if $f \in
\mathcal{F}(U)$, is such that $[f,g]_{U \cap V}=0 $, for any $g
\in
\mathcal{F}(V) $ and any open set of $V$, then $f $ is constant
in the connected components of $U$.
\end{definition}

\begin{example}
\label{sheaf standard manifold}
\normalfont Any  Leibniz manifold $(P, [ \cdot , \cdot ])$ has a
natural presheaf of Leibniz algebras on its presheaf 
of smooth functions that associates to any open subset $U $ of
$P$ the restriction $[ \cdot , \cdot ]| _U $ of $[ \cdot , \cdot
] $ to $C^{\infty} (U)
\times C^{\infty} (U) $.
\end{example}

\subsection{Leibniz reduction by pseudogroups}
\label{Leibniz reduction by pseudogroups subsection}

The main goal of this section is the presentation of a result
that fully characterizes the situations in which the presheaf of
Leibniz algebras in Example~\ref{sheaf standard manifold}
behaves properly under restriction to subsets and projection to
the orbit spaces of pseudogroups of local Leibniz diffeomorphisms
of $(P, [ \cdot , \cdot ])$.

We start by introducing our terminology. Let $P $ be a
smooth manifold and ${\rm Diff} _L (P) $  the
pseudogroup of local diffeomorphisms of $P$. 
More explicitly, the elements of
${\rm Diff} _L (P) $ are diffeomorphisms $F: {\rm
Dom}(F)\subset P
\rightarrow F ({\rm Dom}(F))$ of an open subset ${\rm Dom}(F)
\subset P
$ onto its image $F({\rm Dom}(F))\subset P$. We will denote the
elements of ${\rm Diff} _L (P) $ as pairs $(F, {\rm Dom}(F))$. The
local diffeomorphisms can be composed using the binary operation
defined as
\begin{equation}
\label{composition diffeomorphisms monoid}
\! (G, {\rm Dom}(G))\cdot(F, {\rm Dom}(F)):=(G \circ F, F ^{-1}({\rm
Dom}(G))\cap {\rm Dom}(F)), 
\end{equation}
for all $(G, {\rm
Dom}(G)),(F, {\rm Dom}(F))\in{\rm Diff}_L(P)$. It is easy to see
that this operation is associative and has $(\Bbb I,P)$, the identity
map of
$P$, as (unique) two sided identity element, which makes
${\rm Diff}_L(P)$ into a monoid 
\index{monoid}%
(set with an associative operation
which contains a two--sided identity element). Notice that only the elements
in ${\rm Diff}(P)\subset{\rm Diff}_L(P)$ have an inverse
since, in general, for any $(F, {\rm Dom}(F))\in{\rm Diff}_L(P)$, we
have that 
\begin{eqnarray}
\!\!\!\!\!\!\!\!\!\!\!\!\!\!\!\!\!\!\!\!\!\!\!\!\!\!\!(F^{-1}, F({\rm
Dom}(F)))
\cdot(F,
{\rm Dom}(F))&=&(\Bbb I|_{{\rm Dom}(F)},{\rm Dom}(F))\label{quasigroup 1}\\
\!\!\!\!\!\!\!\!\!\!\!\!\!\!\!\!\!(F, {\rm Dom}(F))\cdot
(F^{-1}, F({\rm Dom}(F)))&=&(\Bbb I|_{F({\rm Dom}(F))},F({\rm
Dom}(F))).\label{quasigroup 2}
\end{eqnarray} 
Consequently, the only way to obtain the identity element $(\Bbb
I,P)$ out of the composition of $F $ with its inverse is having ${\rm
Dom}(F)= P $. It follows from this argument that ${\rm Diff}(P)$ is
the biggest subgroup contained in the monoid ${\rm Diff}_L(P)$ with
respect to the composition law~(\ref{composition diffeomorphisms
monoid}). 
In the sequel we will frequently encounter submonoids
$A$ of
${\rm Diff}_L(P)$ that
satisfy the following property: 
\begin{description}
\item [(PS)]for any
$F: {\rm Dom}(F)\rightarrow  F({\rm Dom}(F))$ in $A$ there exists
another element $F^{-1}:
F({\rm Dom}(F))\rightarrow  {\rm Dom}(F)$ also in $A$ that
satisfies the identities~(\ref{quasigroup 1}) and~(\ref{quasigroup
2}). 
\end{description}
Such submonoids will be referred to as {\bfi pseudogroups}
of  ${\rm Diff}_L(P)$. Recall that $A $ being a submonoid implies that it 
is closed under
composition and
$(\Bbb I, P)
\in A$. 
One of the important features of pseudogroups is that they have
an associated  orbit space. Indeed, if
$A$ is a pseudogroup we define the {\bfi  orbit}
$A \cdot m $ under $A$ of any
element $m \in P $  as the set $A \cdot m:= \{F (m)\mid F
\in A,\, \text{such~ that} \; m \in {\rm Dom}(F)\} $. $A$ being a pseudogroup
implies that the relation
\emph{being in the same $A$--orbit} is an equivalence relation and induces 
a partition
of $P$ into $A$--orbits. The {\bfi  space of $A$--orbits} 
\index{orbit!space}%
\index{space!of orbits}%
will be denoted by $P/A $.
If we endow the space of orbits $P/A $ with the quotient topology, 
the projection $\pi_A:P
\rightarrow P /A$ is a continuous 
and
open map.

Let $S \subset  P $ 
be a subset of $P$ endowed with a topology ${\cal T} $ that in
general does not coincide with the relative or subspace topology.
The presheaf $C^{\infty} _P $  of smooth functions on $P$ induces
a quotient presheaf  $C ^\infty _{P/A }$ on the orbit space $P/A
$. Consider now the subset
\[ A _S:=\left\{ a \in A\mid a (s) \in S\text{ for any } s \in
S\cap {\rm Dom}(a)\right\}.
\] 
All along this section we will assume that $A _S $ is a
subpseudogroup of $A$. This hypothesis will allow us  to
construct the quotients $S/A _S
$ and $P/ A _S $. Given that the quotient $S/A _S $ can be seen
as a subset of $P/A _S $, there is a well defined presheaf of
Whitney smooth functions $ W^{\infty}_{S/ A _S}$ on $S/A
_S $ induced by $C^{\infty} _{P/ A _S} $. We recall
(see~\cite{hsr})  that for any open set $V \subset S/ A _S$, the
elements $f  \in W^{\infty}_{S/ A _S}  (V) $ are characterized by
the fact that if
$\pi_S:S
\rightarrow S/ A _S $ is the projection onto orbit space then for
any $m \in \pi_S ^{-1} (V) $ there exists an open  $A _S
$--invariant  neighborhood of  $m$ in $P$ and $F \in  C^{\infty}
_P (U _m)^{A _S} $ such that 
\begin{equation}
\label{characterization local extension} f \circ \pi _S|_{\pi _S
^{-1} (V)\cap U _m}=F|_{\pi _S ^{-1} (V)\cap U _m}.
\end{equation} We will say that $F$ is a {\bfi  local extension} 
of $f \circ  \pi_S $ at the point $m$. 

\begin{definition} Let $P $ be a smooth manifold, $A \subset {\rm
Diff} _L (P) $  a pseudogroup of local diffeomorphisms of $P$,
and $S$  a subset of $P$ endowed with a topology ${\cal T}$
that is stronger than the relative topology. We say that the
presheaf $W^{\infty}_{S/A _S}$ has the $(A, A _S)$--{\bfi  local
extension property} 
when $A _S$ is a subpseudogroup
of $A$ and for any $f \in W^{\infty}_{S/A _S}(V) $ and $m 
\in \pi_S ^{-1} (V) $ there exists an open $A$--invariant
neighborhood $U _m $ of $m$  in
$M$ and $F \in C^{\infty} _P(U _m) ^A $ such that 
\[ f \circ \pi _S|_{\pi _S ^{-1} (V)\cap U _m}=F|_{\pi _S ^{-1}
(V)\cap U _m}.
\] 
We will say that $F$ is a $A$--{\bfi  invariant local
extension} 
of $f \circ  \pi_S  $ at $m$.
\end{definition}

\begin{definition}
\label{Leibniz reducible definition} 
Let $(P, [ \cdot , \cdot ])$
be a smooth Leibniz manifold and $A \subset \mathcal{P}_L (P)
$  a pseudogroup of local diffeomorphisms of $P$ such that the
presheaf of $A$--invariant functions on $P$ is closed under the
Leibniz bracket $[ \cdot , \cdot ] $. Let
$S \subset P $ be a subset of $P$  such that $ W^{\infty}_{S/A
_S} $ has the $(A, A _S)$--local extension property. We say that
$(P, [ \cdot , \cdot ], A, S)$ is {\bfi  Leibniz reducible} when
$\left(S/A _S, W^{\infty}_{S/A _S}, [ \cdot , \cdot ]^{S/A _S}
\right)$ is a well defined presheaf of Leibniz algebras where,
for any open set $V \subset S/A _S$, the bracket $[ \cdot , \cdot
]^{S/ A _S}_V:W^{\infty}_{S/A _S} (V) \times W^{\infty}_{S/A _S}
(V) \rightarrow W^{\infty}_{S/A _S} (V)  $ is given by 
\begin{equation}
\label{definition Leibniz reduced bracket} [f,g] _V^{S/A
_S}(\pi_S (m))=[F,G] (m)
\end{equation} for any $m \in \pi_S ^{-1} (V) $ and where $F, G $
are $A$--invariant local extensions at
$m$ of
$f \circ \pi_S $ and $g \circ  \pi_S $, respectively.
\end{definition}

The following theorem generalizes the main reduction result
in~\cite{poisson reduction singular} to the context of Leibniz
manifolds with locally defined Leibniz weak symmetries.

\begin{theorem}
\label{main theorem Leibniz reduction} Let $(P, [ \cdot , \cdot
])$ be a smooth Leibniz manifold and $A \subset \mathcal{P}_L (P)
$  a pseudogroup of local diffeomorphisms of $P$ such that the
presheaf of $A$--invariant functions on $P$ is closed under the
Leibniz bracket $[ \cdot , \cdot ] $. Let  $S \subset P
$ be a subset of $P$  such that $ W^{\infty}_{S/A _S} $ has the
$(A, A _S)$--local extension property. Let $B^{\sharp}_L ,
B_R^{\sharp}: T ^\ast M \rightarrow  TM $ be the left and right 
bundle maps respectively  associated to the Leibniz tensor of
$(P,[ \cdot , \cdot ])$. Then
$(P,[\cdot ,\cdot ], A, S)$ is Leibniz reducible if and only if
for any $m \in S $ we have that 
\begin{equation}
\label{fundamental primitive reduction relation}
B_L^{\sharp}(\Delta _m) + B_R^{\sharp}(\Delta _m)\subset 
\left[\Delta _m ^S\right]^{\circ},
\end{equation} where $ \Delta _m:=\{ \mathbf{d} F (m)\mid F \in
C^{\infty} _P(U _m)^A$,  for  any open $A$--invariant
neighborhood $U _m $ of $m$ in $P\}$, and where $\Delta_m ^S=\{
\mathbf{d} F (m) \in \Delta _m\mid F|_{U _m\cap V _m}\text{ is
constant}$, for an open
$A$--invariant neighborhood $U _m $ of  $m$ in $P$  and an open
$A _S $--invariant neighborhood $V _m $ of $m$  in $S\}$.
\end{theorem}

\begin{remark}
\normalfont 
If $S$ has the relative topology then $\Delta_m ^S=\{
\mathbf{d} F (m) \in \Delta _m\mid F|_{U _m\cap S}\text{ is
constant}\,\} $, for an open
$A$--invariant neighborhood $U _m $ of  $m$ in $P$.
\end{remark}

\begin{remark}
\normalfont
If $A$  consists of local Leibniz diffeomorphisms then the
condition on the
presheaf of $A$--invariant functions on $P$ being closed under the
Leibniz bracket $[ \cdot , \cdot ] $ is automatically satisfied.
\end{remark}

\begin{lemma}
\label{prepare for projection function} Let $P$ be a smooth
manifold, $A \subset \mbox{\rm Diff}_L (P) $  a pseudogroup of
local transformations of $P$, and $S\subset P$  a subset whose
topology is stronger than the relative topology and such that $A
_S $ is a subpseudogroup of $A$. If $\pi_S:S
\rightarrow S/A _S $ is the projection, $U \subset P $ is an open
$A$--invariant subset of $P$, $F \in C^{\infty} _P (U) ^A  $, and
$V:= \pi_S (U\cap S)$ then there exists a unique function $f \in
W^{\infty}_{S/A _S}(V) $ such that
\begin{equation}
\label{you will see no invariant} f \circ \pi_S|_{U\cap S}= f
\circ  \pi_S|_{ \pi_S ^{-1} (V)\cap U}=F|_{\pi_S ^{-1} (V)\cap U}.
\end{equation}
\end{lemma}

\noindent\textbf{Proof.\ \ } Since by hypothesis the topology of
$S$ is stronger than the relative topology we have that for any
open $A$--invariant subset $U$ of  $P$, the intersection $U\cap S
$ is an open $A _S$--invariant subset of $S$. As the projection
$\pi_S $ is an open map, the set $V:= \pi_S (U\cap S)$ is open in
$S/A _S $. Also, the $A _S$--invariance of $U \cap S $ implies
that $U\cap S= \pi_S ^{-1} (V)$ and hence
\begin{equation}
\label{you will see no invariant Lemma}
\pi_S ^{-1} (V)\cap U=U\cap S\cap U= U\cap S,
\end{equation} which proves the first equality in~(\ref{you will
see no invariant}).

Now, the invariance properties of  $F$ and $S$ imply the
existence of a unique map $f$ defined on $V$ such that $f \circ 
\pi_S|_{U\cap S}=F|_{U\cap S}$ or equivalently, by~(\ref{you will
see no invariant Lemma}),  $f \circ  \pi_S|_{\pi_S ^{-1} (V)\cap
U}=F|_{\pi_S ^{-1} (V)\cap U}$. Given that by construction $\pi_S
^{-1}(V) \subset U $ then for any $m \in
\pi_S ^{-1} (V) $, the map $f$ satisfies~(\ref{characterization
local extension}) by taking in that characterization $U$ and 
$F$, which implies that $f \in W^{\infty}_{S/A _S} (V)
$.
\quad
$\blacksquare$

\medskip

\noindent\textbf{Proof of Theorem~\ref{main theorem Leibniz
reduction}.\ \ } We first show that if $(P,[ \cdot , \cdot ],
A,S)$ is Leibniz reducible then $ \Delta_m ^S
\subset  \left[ B_L ^{\sharp}(\Delta _m)+ B_R ^{\sharp}(\Delta
_m)\right] ^{\circ}$, for all $ m \in S $. Let
$\alpha_m \in \Delta _m^S $; by definition there exists an open
$A$--invariant neighborhood $U _m $  of $m$ in $P$ and a function
$K \in C^{\infty}_P (U _m)^A $ such that $\alpha_m= \mathbf{d} K 
(m) $ and $K|_{V _m\cap U _m}$ is constant for an open $A _S
$--invariant  neighborhood of $m$ in $S$.  Notice now that by
definition any element in $B_L ^{\sharp}(\Delta_m)+B_R
^{\sharp}(\Delta_m)$ can be written as $X^L _F(m) +X^R_G(m)$ with
$F,G \in C^{\infty}(W _m)^A$,
$W _m $ an open $A$--invariant neighborhood of $m$ in $P$. By
Lemma~\ref{prepare for projection function} there exist functions
$k \in W^{\infty}_{S/A _S}(\pi_S(U _m\cap S))$ and  $f \in
W^{\infty}_{S/A _S}(\pi_S(W _m\cap S))$ such that 
\[ k \circ  \pi_S|_{U _m\cap S}=K|_{U _m\cap S}, \qquad f \circ 
\pi_S|_{W _m\cap S}=F|_{W _m\cap S}.
\] Hence, by the Leibniz reducibility of $(P, [ \cdot , \cdot ],
A,S) $ we have that 
\[
\langle\alpha_m, X^L _F (m) + X^R_G(m)\rangle=[K,G] (m)-
[F,K](m)=[ k| _W, g| _W] _W^{S/A _S}(\pi_S (m))-[ f| _W, k| _W]
_W^{S/A _S}(\pi_S (m)),
\] where $W= \pi_S(U _m\cap S)\cap \pi_S(W _m\cap S)$. However,
given that the function $C$ on $P$ that is constant and equal to
$K (m)  $ is also an
$A$--invariant local extension of $k \circ \pi_S $ at $m$, we
have that 
\[ [ k| _W, g| _W] _W^{S/A _S}(\pi_S (m))-[ f| _W, k| _W] _W^{S/A
_S}(\pi_S (m))=[C,G] - [F,C]=0,
\] which implies that $\langle \alpha_m, X^L _F
(m)+X^R_G(m)\rangle =0 $. Since $X^L _F (m) +X^R_G(m)\in B_L
^{\sharp} (\Delta_m) + B_R^{\sharp} (\Delta_m)$ is arbitrary we
have that $\alpha_m \in \left[ B_L ^{\sharp}(\Delta _m)+
B_R^{\sharp} (\Delta_m)\right] ^{\circ}$.

Suppose now that the inclusion~(\ref{fundamental primitive
reduction relation}) holds and then we will prove the
reducibility of $(P,[\cdot ,\cdot ], A, S)$. Let $f,g \in
W^{\infty}_{S/A _S} (V) $ and $F,G \in C^{\infty} _P (U _m)^A $
be local $A$--invariant extensions of $f
\circ  \pi_S $ and $g \circ \pi_S $, respectively, at a point $m
\in  \pi_S ^{-1} (V) $. We now show that the equality
\begin{equation}
\label{we bracket} [f,g] _V^{S/A _S}(\pi_S (m))=[F,G]_{U _m} (m)
\end{equation} provides a well defined presheaf of Leibniz
algebras. The only point that requires a proof is that the
expression~(\ref{we bracket})  does not depend on the local
extensions utilized in the definition. The fact that $[\cdot ,
\cdot ] ^{S/A _S}$ determines a presheaf of Leibniz algebras is
inherited from the properties of the bracket $[ \cdot , \cdot ] $
on $P$. Let $G ' \in  C^{\infty}_P(U _m) ^A$ be another local
extension of
$g \circ \pi_S $ at $m$. This implies that $G-G'|_{ \pi_S ^{-1}
(V)\cap U _m}=0$ and hence $\mathbf{d}(G-G')(m) \in  \Delta _m ^S
\subset [B_L^{\sharp}(\Delta _m)+B_R^{\sharp}(\Delta
_m)]^{\circ}$. Consequently,
\[ 0=\langle \mathbf{d}(G-G') (m) , X _F ^L
(m)\rangle=-[F,G-G']_{U _m} (m),
\] which implies that $[F,G]_{U _m} (m)=[F, G']_{U _m} (m)$ and
hence guarantees the independence of~(\ref{we bracket}) with
respect to the choice of local extension for $g \circ \pi_S $. A
similar argument guarantees that this definition is also
independent of the choice of extension for $f \circ  \pi_S$.
Therefore, the expression~(\ref{we bracket}) defines a function
$[f,g] _V^{S/A _S} $ on
$V$ that actually belongs to $W^{\infty} _{S/A _S} (V) $ because
if $F $ and $G$ are local $A$--invariant extensions of $f
\circ  \pi_S $ and $g \circ \pi_S $, respectively, at any point
$m \in  \pi_S ^{-1} (V)
$ then 
so is the function $\{F,G\} $ with respect to $\{f,g \}^{S/A} _V
\circ  \pi_S$ by the
hypothesis on the
presheaf of $A$--invariant functions on $P$ being closed under the
Leibniz bracket $[ \cdot , \cdot ] $.
\quad $\blacksquare$

\subsection{Leibniz reduction by distributions}

The Leibniz reduction theorem that we presented in
Section~\ref{Leibniz reduction by pseudogroups subsection}
requires the presence of a pseudogroup of transformations defined
in the entire manifold. However, sometimes one may want to reduce
with respect to an invariance property defined only on a subset
of the manifold in question. The study of this situation is the
main goal of this section. We start by introducing the setup.

\begin{definition} Let $P$ be a differentiable manifold and
$S\subset P$  a  decomposed subset of $P$. Let 
$\{S_i\}_{i\in I}$ be the pieces of this decomposition. The
topology of $S$ is not necessarily the relative topology as a
subset of $P$.  We say that $D\subset TP|_S$ is a  {\bfi  smooth
distribution on $S$ adapted to the decomposition\/} 
$\{S_i\}_{i\in I}$, if
$D\cap TS_i$ is a smooth distribution on $S_i$ for all $i\in I$.
The distribution $D$ is said to be {\bfi  integrable\/}  if
$D\cap TS_i$ is integrable for each $i\in I$.
\end{definition}

In the situation described by the previous definition and if $D$
is integrable,  the  integrability of the distributions $D_{S
_i}:=D\cap TS_i$ on $S_i$ allows us to partition each $S_i$ into
the corresponding maximal integral manifolds.  Thus, there is an
equivalence relation  on $S_i$ whose  equivalence classes are
precisely these maximal integral manifolds. Doing this on each
$S_i$, we obtain an equivalence  relation $D_S$ on the whole set
$S$ by taking the union of the different  equivalence classes
corresponding to all the $D_{S _i}$. We define the quotient space
$S/D _S$ as
\[S/D _S:=\bigcup_{i\in I} S_i/D_{S _i}.\]   We will denote by
$\pi_{D_S }: S \rightarrow  S/ D_S $ the natural projection.

\begin{definition}
\label{Leibniz distribution weak} Let $(P,\,[\cdot ,\cdot])$ be a
Leibniz manifold and 
$D\subset TP$  a smooth distribution on $P$. The distribution
$D$ is called {\bfi  Leibniz\/} or {\bfi  canonical\/}, if the
condition $\mathbf{d} f|_D=\mathbf{d} g|_D=0$, for any $f,\,g\in
C_P^{\infty}(U)$ and any open subset $U \subset P $, implies that
$\mathbf{d} [f,\,g]|_D=0$. 
\end{definition}

\bigskip

\noindent {\bf The presheaf of smooth functions on $S/D_S$.} In
this section we will be considering a presheaf of smooth
functions on $S/ D_S $  that require less invariance properties
in their definition than those that appeared in the context of
quotients by pseudogroups of transformations. We define the
presheaf of smooth functions
$C^{\infty} _{S/ D_S} $ on $S/ D _S $ as the map that associates
to any open subset $V $ of $S/ D _S$ the set of functions
$C^{\infty} _{S/ D_S} (V)$ characterized by the following
property: $f  \in C^{\infty} _{S/ D_S} (V)$ if and only if for
any $z \in V $  there exists $m \in \pi _{D_S} ^{-1} (V) $, $U _m
$ open neighborhood of $m$ in  $P$, and
$F \in C^{\infty} _P (U _m) $ such that 
\begin{equation}
\label{one more presheaf} f \circ  \pi_{D _S}|_{\pi _{D_S} ^{-1}
(V)\cap U _m}=F|_{\pi _{D_S} ^{-1} (V)\cap U _m}.
\end{equation} We say that $F$ is a {\bfi  local extension} of $f
\circ \pi_{D_S} $ at the point $m \in \pi_{D_S} ^{-1} (V) $. It
can be proved (see~\cite{hsr}) that if $S$ is a smooth embedded
submanifold of $P$ and
$D_S $ is a smooth, integrable, and regular distribution on $S$
then the presheaf
$C^{\infty}_{S/ D_S}$ coincides with the presheaf of smooth
functions on $S / D_S $  when considered as a regular quotient
manifold.

We say that the presheaf $C^{\infty}_{S/ D_S}$ has the $(D,
D_S)$--{\bfi  local extension property} when the topology of $S$
is stronger than the relative topology and, at the same time, the
local extensions of
$f 
\circ
\pi_{D_S}
$ defined in~(\ref{one more presheaf}) can always be chosen so
that
\[
\mathbf{d} F (n) |_{D(n)}=0, \quad \mbox{for any}\quad n \in  \pi
_{D_S} ^{-1} (V)\cap U _m.
\] 
We say that $F$ is a {\bfi  local $D$--invariant extension}
of $f \circ \pi_{D_S} $ at the point $m \in \pi_{D_S} ^{-1} (V) $.

\begin{definition} Let  $(P, [ \cdot , \cdot ])$ be a Leibniz
manifold, $S$  a decomposed subset of $P$, and $D \subset  TP|
_S $  a Leibniz integrable generalized distribution adapted to
the decomposition of $S$. Assume that  $C^{\infty} _{S/ D_S}$ has
the $(D, D_S)$--local extension property. We say that $(P, [
\cdot , \cdot ], D, S)$ is {\bfi  Leibniz reducible}  when  $(S/
D_S, C^{\infty}_{S/ D_S},[ \cdot , \cdot ]^{S/ D_S})$ is a well
defined presheaf of Leibniz algebras where, for any open set $V
\subset S/ D_S$, the bracket $[
\cdot , \cdot ]_V ^{S/ D_S}: C^{\infty} _{S/ D_S} (V) \times
C^{\infty} _{S/ D_S} (V)
\rightarrow C^{\infty} _{S/ D_S} (V) $ is given by
\[ [ f , g ]_V ^{S/ D_S}(\pi_{D_S} (m)):=[F, G] (m),
\] for any $m \in \pi_{D_S}^{-1} (V) $. The maps $F,G $ are local
$D$--invariant extensions at $m$ of $f \circ \pi_{D_S} $ and $g
\circ \pi_{D_S} $, respectively.
\end{definition}

The proof of the following theorem mimics the corresponding
implication in theorem~\ref{main theorem Leibniz reduction}.

\begin{theorem}
\label{reduction distribution singular only one} Let $(P, [ \cdot
, \cdot ]) $ be a Leibniz manifold with associated Leibniz tensor
$B$ , $S$  a decomposed space, and $D \subset  TP| _S
$  a Leibniz integrable generalized distribution adapted to the
decomposition of $S$. Assume that  $C^{\infty} _{S/ D_S}$ has the
$(D, D_S)$--local extension property. Then
$(P, [
\cdot , \cdot ], D, S)$ is Leibniz reducible if for any
$m \in S  $
\begin{equation}
\label{Leibniz reducibility distribution singular}
B_L^{\sharp}(\Delta _m)+B_R^{\sharp}(\Delta _m) \subset  \left[
\Delta _m ^S  \right] ^{\circ}
\end{equation} where $\Delta_m:=\{ \mathbf{d} F (m)\mid F \in
C^{\infty}_P (U _m), \mathbf{d}F (z)|_{D (z)}=0, \text{ for all
}z \in U _m\cap S,$
 and for any open neighborhood 
$U_m  $ of $m$ in $P\}$ and $\Delta _m^S:=\{ \mathbf{d} F (m) \in
\Delta _m\mid F|_{U _m\cap V _m}$ is constant for an open
neighborhood  $U _m    \text{ of }   m  \text{ in }  P $ and an
open neighborhood
$V _m$ of $m$ in $S\}
$.
\end{theorem} 

\begin{remark}
\normalfont If $S$ is endowed with the relative topology then
$$\Delta _m^S:=\{ \mathbf{d} F (m) \in \Delta _m\mid F|_{U _m\cap
V _m}\text{ is constant for an open neighborhood } U _m   
\text{ of }   m  \text{ in }  P\}.
$$
\end{remark}

\begin{remark}
\normalfont As opposed to the situation in Theorem~\ref{main
theorem Leibniz reduction}, the condition~(\ref{Leibniz
reducibility distribution singular}) is sufficient for Leibniz
reducibility but in general not necessary. The reason behind this
circumstance is that the functions that define the spaces $\Delta
_m  $ and $\Delta_m ^S $ are not defined on saturated open sets
which prevents the formulation of a result similar to 
Lemma~\ref{prepare for projection function}. As we will see in
Theorem~\ref{reduction theorem by distributions regular case}, 
an alternative hypothesis that makes this condition necessary and
sufficient is, roughly speaking, the regularity of the
distribution
$D _S:= D\cap TS $.
\end{remark}

\noindent {\bf Reduction by regular canonical distributions.} Let
$(P,[ \cdot , \cdot ])$ be  a Leibniz manifold and $S$   an
embedded submanifold of $P$. Let $D\subset TP| _S $  be a
subbundle of the tangent bundle of $P$ restricted to $S$ such
that 
$D_S:=D\cap TS  $  is a smooth, integrable, and
regular distribution  on
$S$ and 
$D$ is Leibniz.
Our next theorem is a generalization of the main result
of~\cite{poisson reduction marsden ratiu} to the context of
Leibniz manifolds.

\begin{theorem}
\label{reduction theorem by distributions regular case}  Let
$(P,[ \cdot , \cdot ])$ be a Leibniz manifold with associated
Leibniz tensor $B$ and
$S$   an embedded smooth submanifold of  $P$. Let $D\subset 
TP| _S$ be a canonical subbundle of the tangent bundle of $P$
restricted to $S$ such that  $D_S:= D\cap TS $ is a smooth,
integrable, and regular distribution on $S$. Then $(P, [ \cdot ,
\cdot ], D, S) $ is Leibniz reducible if and only if 
\begin{equation}
\label{classical reduction condition Marsden Ratiu} B
_L^{\sharp}(D^{\circ})+B_R ^{\sharp}(D^{\circ})\subset TS + D.
\end{equation}
\end{theorem}

\noindent\textbf{Proof.\ \ } We first prove that the
condition~(\ref{classical reduction condition Marsden Ratiu})
implies the Leibniz reducibility of
$(P, [ \cdot , \cdot ], D, S)
$. This implication can be obtained as a corollary  of
Theorem~\ref{reduction distribution singular only one}.
Indeed, a result whose proof can be found in~\cite{hsr}
guarantees that the
hypotheses  on $D _S $ imply that the presheaf $C^{\infty}_{S/
D_S}
$ has the $(D, D_S)$--local extension property. Hence, it suffices
to show that in this situation 
\begin{eqnarray}
\Delta _m&=& D (m) ^{\circ}, \label{first equality to be shown
13}\\
\left[\Delta _m ^S\right]^{\circ}&=&T _mS+ D(m). \label{first
equality to be shown 14}
\end{eqnarray} In order to prove~(\ref{first equality to be shown
13}) notice first that by definition
$\Delta_m \subset  D (m) ^{\circ} $. To prove the converse
inclusion take $\alpha _m \in  D(m)^{\circ} $ arbitrary and let 
$U
_m
$ be a submanifold chart of $S$  around $m$ that we can think
of as
$U
\times V
\subset 
\mathbf{F}_1\oplus
\mathbf{F} _2 $, where $U $ and $V$ are open neighborhoods of the
origin in two vector spaces $\mathbf{F} _1 $ and $\mathbf{F}_2 $,
respectively. This chart can be constructed so that
$m\equiv(0,0)$ and $U _m\cap S= U
$. Additionally, we can locally take
$D=U
\times  \mathbf{E}$, with $\mathbf{E}  $ a vector subspace of
$\mathbf{F} _1\oplus
\mathbf{F}  _2 $,  $T U _m= U \times  V \times ( \mathbf{F}
_1\oplus \mathbf{F} _2) $, and $T ^\ast  U _m= U \times  V \times
(\mathbf{F} _1\oplus \mathbf{F} _2) ^\ast $. In these coordinates
 $\alpha _m\equiv (0,0, \alpha)$, with $\alpha
\in
\mathbf{E} ^{\circ} $. Define $F : U \times  V \rightarrow 
\mathbb{R} $ by $F(u,v)=
\langle \alpha, (u,v) \rangle $. Note that $\mathbf{d} F (m)
\equiv \mathbf{d}F (0,0)=
\alpha \equiv \alpha _m $ and that for any $(u,0,w) \in D (u,0)
$, $u \in U $, $w \in
\mathbf{E} $ we have that $\mathbf{d} F (u,0 ) \cdot  w= \langle
\alpha, w \rangle =0 $, which implies that $\mathbf{d} F(z)|_{D
(z)}=0 $, for any $z \in U _m\cap S $, as required.

We now prove~(\ref{first equality to be shown 14}). By definition
$ T _mS + D (m)
\subset  [ \Delta _m ^S]^{\circ}$. Conversely, the inclusion  $[
\Delta _m ^S]^{\circ}
\subset T _mS + D (m) $ holds if and only if $ D (m) ^{\circ}\cap
T _mS ^{\circ} \subset 
\Delta _m ^S $ which, by~(\ref{first equality to be shown 13}),
amounts to $\Delta _m
\cap T _mS ^{\circ} \subset  \Delta _m^S $. We prove this
inclusion by using again the same adapted local submanifold
coordinates around the point $m$. Let $\alpha _m \in
\Delta _m
\cap T _mS ^{\circ} $ arbitrary. As we saw above, there exists
$\alpha \in  \mathbf{E} ^{\circ} $ such that $\alpha _m  =
\mathbf{d} F (0,0)$, with $F: U \times  V \rightarrow 
\mathbb{R}  $ given by $F (u,v):= \langle  \alpha, (u,v) \rangle 
$, $(u,v) \in U \times V $. Since $\alpha _m \in (T _mS) ^{\circ}
$ we have that $\mathbf{d} F (0,0) \cdot (u,0)= \langle  \alpha,
(u,0) \rangle =0 $, for any $u \in \mathbf{F} _1 $. This equality
implies that $F (u,v)= \langle  \alpha, (u,0)\rangle + \langle 
\alpha, (0,v)\rangle =
\langle  \alpha, (0,v)\rangle  $, for any
$(u,v)
\in U
\times V $, and hence $F $ is constant on $U=U _m\cap S $  which
shows that $\alpha _m
\in \Delta _m^S $, as required.

We now show that the reducibility of  $(P, [ \cdot , \cdot ], D,
S)$ implies~(\ref{classical reduction condition Marsden Ratiu})
or, equivalently, that for any $m \in S $
\begin{equation}
\label{another way to see the thing} T _mS ^{\circ}\cap D
(m)^{\circ} \subset  \left( (B_L ^{\sharp } (m)+B_R ^{\sharp }
(m))  (D (m) ^{\circ})
\right)^{\circ}. 
\end{equation} We proceed again by using the same local
coordinates. In this occasion we consider a non degenerate inner
product $\langle  \cdot , \cdot  \rangle _{\mathbf{F} _1\oplus
\mathbf{F} _2} $ on $\mathbf{F} _1\oplus \mathbf{F} _2 $ defined
by $\langle  (u _1, u _2) , (v _1, v _2)  \rangle _{\mathbf{F}
_1\oplus \mathbf{F} _2}= \langle u _1, v _1 \rangle_{ \mathbf{F}
_1}+
\langle  u _2, v _2 \rangle  _{\mathbf{F}_2}  $, with $\langle 
\cdot , \cdot  
\rangle _{\mathbf{F} _1} $ and $\langle  \cdot , \cdot  \rangle
_{ \mathbf{F} _2} $ non degenerate inner products in $\mathbf{F}
_1 $ and $\mathbf{F} _2 $, respectively, and $u _1, v _1 \in
\mathbf{F} _1 $, $u _2, v _2 \in  \mathbf{F} _2 $. Given that $U
_m \cap S= U\subset
\mathbf{F}  _1$, any element $\alpha _m \in T _mS ^{\circ} $ can
be written as $ \alpha _m= \langle (0, u _0), \cdot  \rangle
_{\mathbf{F} _1\oplus \mathbf{F} _2} $, for some
$u _0 \in  \mathbf{F} _2$ or, analogously, as $\alpha _m=
\mathbf{d}K (m) $, with $K \in C^{\infty} (U _m)$ defined by 
\begin{equation}
\label{definition of k k k} K(u,v):= \langle (0, u _0), (u,v)
\rangle _{\mathbf{F} _1\oplus \mathbf{F} _2}= \langle  u _0, v
\rangle _{\mathbf{F} _2}.
\end{equation}  Moreover, if $\alpha_m \in T _mS ^{\circ}\cap D
(m) ^{\circ} $ then as $D$ in these coordinates looks like $D=U
\times \mathbf{E} $ for some vector subspace $\mathbf{E}
\subset  \mathbf{F} _1\oplus \mathbf{F} _2 $, we have that the
function $K $ defined in~(\ref{definition of k k k}) is such that 
\[ K|_{U _m\cap S}=0\text{ and } \mathbf{d}K (z)|_{D
(z)}=0,\text{ for any } z \in U _m\cap S.
\] We have thus proven that any $\alpha _m \in T _m S
^{\circ}\cap D (m) ^{\circ} $ can be written as $\alpha _m=
\mathbf{d}K (m)  $ with $K$  a local $D$--invariant extension of
the zero function in $S$ at the point $m \in S $.

Let now $\beta_m \in  D (m) ^{\circ} $. Due to the non degeneracy
of the inner product $\langle  \cdot , \cdot  \rangle
_{\mathbf{F} _1\oplus \mathbf{F} _2} $ there exists $w _0 \in
\mathbf{F} _1\oplus \mathbf{F} _2$ such that  $\beta _m=
\mathbf{d} F (m) $ with $F (u):= \langle w _0, u \rangle
_{\mathbf{F} _1\oplus \mathbf{F} _2}$, $ u \in U _m $, and such
that  $\langle  w _0, w \rangle _{\mathbf{F} _1\oplus \mathbf{F}
_2}=0$ for any $w \in \mathbf{E} $. The regularity of the
distribution $D_S $ implies via a result of Godement (see Lemma
3.5.26 in~\cite{mta}) that the neighborhood $U _m $ can be shrunk
so that there exists a smooth submanifold $T$ of $U _m \cap S $
(called a {\bfi  local slice} of $D _S $) and a smooth map $s: U
_m\cap S \rightarrow  T $ such that  $s|_{T} $ is the identity
map on $T$ and the integral leaf $\mathcal{L} _{u} $ of $D_S $
that contains any arbitrary point $u \in U _m\cap S $ is such
that  $\mathcal{L} _{u}\cap T=\{s (u)\} $. Notice now that since
$\mathbf{d} F|_{U _m\cap S}\cdot  D _S|_{U _m\cap S}=0$, we can
use the slice and the map $F|_{U _m\cap S} $ to define another
map $f \in C^{\infty}_{S/ D_S}(\pi_{D_S}(U _m\cap S)) $ as the
unique map that satisfies
\[ f \circ  \pi_{D_S}(z)=F (z)=F (s (z)),\text{ for any } z \in U
_m\cap S.
\]

Using the constructions in the last two paragraphs we can now
write for any $\alpha _m
\in T _m S ^{\circ}\cap D (m) ^{\circ} $ and any $\beta^1_m,
\beta_m^2 \in  D (m) ^{\circ} $
\begin{multline}
\label{we now kill it}
\langle \alpha _m, (B_L ^{\sharp} (\beta^1_m)
+B_R^{\sharp}(\beta^2 _m) \rangle \\ =[K,F] (m)-[G,K]=[0, f] ^{S/
D_S}_{\pi_{D_S} (U _m\cap S)}(\pi_{D_S}(m))-[g, 0]^{S/
D_S}_{\pi_{D_S} (U _m\cap S)}(\pi_{D_S}(m))\\ =[0,F](m)-[G,0]
(m)=0,
\end{multline} where in the last equality we used the Leibniz
reducibility of $(P, [ \cdot , \cdot ], D, S)$ to write
\[[0,f]_{\pi_{D_S}(U _m\cap S)}^{S/ D_S}(\pi_{D_S} (m))=[0,F]
(m)\quad \text{and} \quad 
[g,0]_{\pi_{D_S}(U _m\cap S)}^{S/ D_S}(\pi_{D_S} (m))=[G,0] (m),
\] 
since the zero function on $M$ is also an extension of the zero
function on $S$ that can be used instead of $K$ in the definition
of the bracket. The expression~(\ref{we now kill it}) 
establishes~(\ref{another way to see the thing}). \quad
$\blacksquare$ 

\bigskip

\addcontentsline{toc}{section}{Acknowledgments}
\noindent\textbf{Acknowledgments.}  We thank Richard Cushman, 
Mark Roberts, and Miguel Rodr\'{\i}guez--Olmos for many
illuminating  discussions. This research
was partially supported by the European Commission through
funding for the Research Training Network
\emph{Mechanics and Symmetry in Europe} (MASIE) and by the Marie Curie 
fellowship
HPTM-CT-2001-00233.

\end{document}